\documentclass{article}

\usepackage{amsmath}
    \numberwithin{equation}{section}
\usepackage{amssymb}
\usepackage{amsthm}
\usepackage{mathtools}

\usepackage{adjustbox}
\usepackage{arydshln}
\usepackage{booktabs}
\usepackage{caption}
\usepackage{float}
\usepackage{pgfplots}
    \pgfplotsset{width=15.3cm,compat=1.9}
\usepackage{subcaption}
\usepackage{tikz}
    \usetikzlibrary{plotmarks}
    \usetikzlibrary{arrows.meta}
    \usetikzlibrary{shadings,intersections}
\usepackage{tikz-qtree}
\usepackage{xcolor}
\usepackage{multirow}
\usepackage{multicol}
\usepackage{wrapfig}
\usepackage{subcaption}

\usepackage[english]{babel}
\usepackage[top=2.5cm, left=3cm, right=3cm, bottom=4cm]{geometry}
\usepackage{graphicx}
    \graphicspath{ {./images/} }
\usepackage[utf8]{inputenc}
\usepackage{parskip}
\usepackage{todonotes}

\usepackage[colorlinks=true, allcolors=black]{hyperref}
\usepackage{cleveref}
\usepackage{cite}
\newcommand{\tref}[1]{Table \ref{t:#1}}

\newcommand{\sref}[1]{Section \ref{s:#1}}
\newcommand{\lref}[1]{Lemma \ref{l:#1}}
\newcommand{\thref}[1]{Theorem \ref{th:#1}} 

\usepackage{blindtext}
\usepackage{ragged2e}

\newcommand{\iZ}{\in\mathbb{Z}}
\newcommand{\N}{\mathbb{N}}

\newcommand{\cc}[1]{|{#1}|}
\newcommand{\ccc}[1]{||{#1}||}
\newcommand{\cccc}[1]{|||{#1}|||}

\newcommand{\spt}{\mathrm{spt}}

\newcommand{\lcm}{\mathrm{lcm}}

\newcommand\sums{%
  \operatornamewithlimits{%
    \mathchoice
     {\vcenter{\hbox{\huge $\xi$}}}
     {\vcenter{\hbox{\Large $\xi$}}}
     {\xi}
     {\xi}}}

\newcommand\TikCircle[1][2.5]{\tikz[baseline=-#1]{\draw(0,0)circle[radius=#1mm];}}
\newcommand\circs{%
  \operatornamewithlimits{%
    \mathchoice
     {\vcenter{\hbox{\huge $\TikCircle[4]$}}}
     {\vcenter{\hbox{\Large $\TikCircle[4]$}}}
     {\TikCircle[4]}
     {\TikCircle[4]}}}
\newcommand\TikSquare[1][2.5]{\tikz[baseline=-#1]{\draw(-0.4,0.4)rectangle(0.4,-0.4);}}
\newcommand\squares{%
  \operatornamewithlimits{%
    \mathchoice
     {\vcenter{\hbox{\huge $\TikSquare[4]$}}}
     {\vcenter{\hbox{\Large $\TikSquare[4]$}}}
     {\TikSquare[4]}
     {\TikSquare[4]}}}

\newtheorem{theorem}{Theorem}[section]
\newtheorem{lemma}[theorem]{Lemma}
\newtheorem{corollary}[theorem]{Corollary}

\newtheorem{conjecture}{Conjecture}[section]
\theoremstyle{remark}
\newtheorem*{remark}{Remark}

\title{Closed-Form Formula for the Partition Function and Related Functions}
\author{Alfredo Nader}
\date{}

\begin{document}

\maketitle

\begin{abstract}
    We develop a new closed-form arithmetic and recursive formula for the partition function and a generalization of Andrews' smallest parts (spt) function. Using the inclusion-exclusion principle, we additionally develop a formula for the not-relatively prime partition function (which counts the number of partitions that are not relatively prime). Moreover, we prove a theorem involving the greatest common divisor of partitions, which allows us to link partitions to prime numbers and lets us derive a formula for the relatively prime function. Lastly, we develop numerous new identities for Jordan's totient function of second order, Euler's totient function, and Dedekind's psi function.
\end{abstract}



\section{Introduction}
Two integer numbers, $a$ and $b$, are relatively prime if they do not share any common factors, i.e., $\gcd(a,b)=1$, and use the notation $a\perp b$ \cite[Section 4.5]{graham}. For any set $A$, we say that $A$ is relatively prime if $\gcd(A)=1$ \cite{nathanson}, and we denote it by $A\vdash$. We write $A\nvdash$ if $A$ is not relatively prime. We denote the cardinality of some set $A$ by $|A|$. If $|A|=1$ and $x\in A$, then $\gcd(A)=x$. 

A partition of a nonnegative integer $n$ is a way of writing $n$ as a sum of positive integers, without considering the order of the summands. Individual summands in partitions are called parts. The notation $\lambda^n_k=\{\lambda_1,\lambda_2,\dots,\lambda_k\}$ resembles a partition of $n$ with $k$ parts, where $\lambda_1,\lambda_2,\dots,\lambda_k$ are the parts of $\lambda_k^n$. We say that $\lambda_k^n$ is a relatively prime partition if $\gcd(\lambda_k^n)=1$. As opposed to usual, we use the convention that $\lambda_1\leq\lambda_2\leq\cdots\leq\lambda_k$. 

Throughout, the set $\N=\{x\iZ:x>0\}$. $p$ runs through the primes and every variable rather than $p$ runs through $\N$, unless otherwise indicated. We write $a\mid b$ to say that $a$ divides $b$ and write $a\nmid b$ to say that $a$ does not divide $b$. We denote the floor function, the ceiling function, and the nearest integer function by $\lfloor x\rfloor$, $\lceil x\rceil$, and $\langle x\rangle$, respectively. We denote the number of divisors of $n$ by $d(n)$.

If $f(n)$ is a function related to partitions, we will write $f(n,k)$ to express the same function only considering partitions of $k$ parts. We denote by $p(n)$ the partition function of $n$, which counts the number of partitions of $n$, we denote by $p_\Psi(n)$ the relatively prime partition function, which counts the number of relatively prime partitions of $n$, and we denote by $\Lambda(n)$ the not-relatively prime partition function of $n$, which counts the number of partitions of $n$ that are \textit{not} relatively prime. We let the set $\overline{\lambda^n}$ contain all partitions of $n$ and the set $\overline{\lambda_k^n}$ contain all partitions of $n$ with $k$ parts. Note that
\begin{align*}
    \left|\overline{\lambda^n}\right|=p(n)\coloneqq\sum_{\lambda_k^n\in\overline{\lambda^n}}1;\quad\left|\overline{\lambda^n_k}\right|=p(n,k)\coloneqq\sum_{\lambda_k^n\in\overline{\lambda_k^n}}1.
\end{align*}

\begin{lemma}\label{l:primesp}
    $\gcd(\lambda_1^n)=n$ and $p(n,1)=1$ for all $n\geq1$.
\end{lemma}
\begin{proof}
    By definition, $\overline{\lambda_1^n}\coloneqq\{n\}$ since $\lambda_1=n$
    and, therefore, we directly see that $\gcd(\lambda_1^n)=n$ and hence $p(n,1)=1$ for any $n\geq1$.
\end{proof}

If $f(n)$ is a function related to partitions, we use the convention that $f(n)=f(n,k)=0$ for any $n<0$, $k<1$, or $k>n$. Hence, note that
\begin{align*}
    f(n)=\sum_{k}f(n,k)=\sum_{k<1}f(n,k)+\sum_{k=1}^nf(n,k)+\sum_{k>n}f(n,k)=\sum_{k=1}^nf(n,k).
\end{align*}

\begin{remark}
    We use the convention $p(0)=p_\Psi(0)=1$ since we consider the (empty) partition $\lambda_1^0=\{\}$ as the only partition of 0. We do this to be able to extend the partition function to all nonnegative integers. However, to be as consistent as possible, we additionally let $p(0,1)=p_\Psi(0,1)=1$.
\end{remark}

Due to \lref{primesp}, 
\begin{equation}\begin{split}\label{eq:base_relprime}
    \Lambda(n,1)&=\begin{cases} 0 & n=0,~n=1 \\ 1 & n>1 \end{cases},\\
    p_\Psi(n,1)&=\begin{cases} 1 & n=0,~n=1 \\ 0 & n>1 \end{cases},\\
    \notag p(n,1)&=p_\Psi(n,1)+\Lambda(n,1)=1
\end{split}\end{equation}
for all $n\geq0$.

We denote Jordan's totient function of order $m$ by
\begin{align*}
    J_m(n)\coloneqq n^m\prod_{p\mid n}\left(1-\frac1{p^m}\right),
\end{align*}
Euler's totient function by
\begin{align*}
    \phi(n)\coloneqq n\prod_{p\mid n}\left(1-\frac1p\right)=J_1(n),
\end{align*}
and Dedekind's psi function by
\begin{align*}
    \psi(n)\coloneqq n\prod_{p\mid n}\left(1+\frac1p\right)=\frac{J_2(n)}{J_1(n)}=\frac{J_2(n)}{\phi(n)},
\end{align*}
where $p$ runs through the prime factors of $n$ in all three cases.

\begin{theorem}\label{th:part}
    For any $1\leq k\leq n$, $\gcd(\lambda_k^n)=m$ if and only if $m\mid n$.
\end{theorem}

In \sref{th} we prove \thref{part}. In \sref{formula} we look at tables containing partitions with $k\leq4$ parts, and we use them to develop a recursive and a closed-form formula for $p(n,k)$ and $p(n)$. In \sref{spt} we generalize Andrew's smallest parts function (see \cite{andrews}) and offer a recursive and closed-form formula for it. Using the inclusion-exclusion principle, we develop a formula for $\Lambda(n,k)$ and $\Lambda(n)$ in \sref{lambda}. Additionally, we use Corollary \ref{c:formula} to develop a formula for $p_\Psi(n,k)$ and $p_\Psi(n)$. Lastly, in \sref{imp} we utilize some connections that have been made with the relatively prime partition function in \cite{elbachraoui} to derive new identities for Jordan's second-order totient function, Euler's totient function, and Dedekind's psi function.

To express the formula we will develop in \sref{formula} for the partition function, we need a new sum notation, which we will call \textit{xi-notation} or \textit{$\xi$-notation}. A sum specified by two or more indices is called a multiple sum and, more generally, a sum that is specified by $n$ indices is called a $n$-multiple sum. Let $P(j)$ be some Boolean property, and let $Q(j)$ be the number of $j$ such that $P(j)=\text{True}$. Similar to how we can write
\begin{align*}
    \sum_{P(j)}a_j&=\overbrace{a_1+\cdots+a_{Q(j)}}^{Q(j)\text{ terms}}
\end{align*}
to express a sum of all terms $a_j$ such that $j$ is an integer satisfying $P(j)$ \cite[p. 23]{graham}, we may write
\begin{align*}
    \sums_{P(j)}S=\overbrace{\sum_{P_1(\nu)}\sum_{P_2(\nu)}\cdots\sum_{P_{Q(j)}(\nu)}}^{Q(j)\text{ sums}}S.
\end{align*}
to express a $Q(j)$-multiple sum in which the index of each sum is some Boolean property $P_i(\nu)$ and the summand of the multiple sum is $S$. We can use brackets below the index of the multiple sum to specify certain parameters that the sums must have, that is, we can write, for example,
\begin{align}
    \notag\sums_{{\substack{a\leq j\leq b \\ [\tau,~\nu_1,~\nu,~\ell]}}}S&=\overbrace{
    \sum_{\nu_1\leq m_\tau\leq\ell}
    \,\sum_{\nu\leq m_\tau\leq\ell}\cdots
    \,\sum_{\nu\leq m_\tau\leq\ell}}^{b-a+1\text{ sums}}S \\
    \label{eq:xinotation}\sums_{{\substack{j=a \\ [\tau,~\nu_1,~\nu,~\ell]}}}^bS&=\overbrace{
    \sum_{m_\tau=\nu_1}^{\ell}
    \sum_{m_\tau=\nu}^{\ell}\cdots
    \sum_{m_\tau=\nu}^{\ell}}^{b-a+1\text{ sums}}S,
\end{align}
where $a\leq b$, $S$ is the summand of the multiple sum, and $\nu$ and $\ell$ are functions for the lower and upper bounds of the sums, respectively. In this sense, $a,~b,~\tau,~\nu_1,~\nu,~\ell,~S$ can depend on any function with any appropriate input to the context in use. This is the notation we will continue using in this work. 

\begin{remark}
    If $a>b$ in \eqref{eq:xinotation} we can say that (1) it is undefined, or (2) it is equal to $S$. Here, we will use the convention of the first option.
\end{remark}

As an example, we may write Creech's formula for the number of $k$-almost primes less than $x$ (see \cite[Theorem 1.1]{creech}), using $\xi$-notation, as
\begin{align*}
    \pi_k(x)&=\sum_{i_{k-1}=1}^{\pi\left(\sqrt[k]x\right)}\sum_{i_{k-2}=i_{k-1}}^{\pi\left(\sqrt[k-1]{\frac{x}{p_{i_{k-1}}}}\right)}\sum_{i_{k-3}=i_{k-2}}^{\pi\left(\sqrt[k-2]{\frac{x}{p_{i_{k-2}}p_{i_{k-1}}}}\right)}\cdots\sum_{i_{1}=i_{2}}^{\pi\left(\sqrt{\frac{x}{p_{i_2}p_{i_3}\cdots p_{i_{k-1}}}}\right)}\left(\pi\left(\frac{x}{p_1p_2\cdots p_{k-1}}\right)-i_1+1\right)\\
    &=\sums_{\substack{j=1 \\ \left[k-j,~1,~i_{k-j+1},~\pi\left(\sqrt[k-j+1]{\prod_{r=k-j+1}^{k-1}\frac{x}{p_{i_r}}}\right)\right]}}^{k-1}\left(\pi\left(\prod_{i=1}^{k-1}\frac{x}{p_{i}}\right)-i_1+1\right).
\end{align*}

\section{\thref{part}}\label{s:th}

\begin{proof}[Proof of \thref{part}]
To prove \thref{part}, we first note from \lref{primesp} that, in general, $\gcd(\lambda_1^n)=n\mid n$, and since 1 is divisible by itself, then $\gcd(\lambda_k^1)=1\mid1$. Subsequently, we directly prove in \lref{primes} that any partition of any prime number with two or more parts is relatively prime. Then, by direct proof and using \lref{primes}, we show in \lref{all} that, in general, for any $m\mid n$ and $1\leq k\leq n$, we \textit{may} have $\gcd(\lambda^n_k)=m$. Lastly, in \lref{iff} we show by contradiction that the only way \lref{all} holds is when $m\mid n$, therefore proving \thref{part}. \thref{part} implies Corollaries \ref{c:prime}, \ref{c:general}, and \ref{c:formula}. Before starting, we need to remark something:

\begin{remark}
    In \lref{all} we prove that for $m\mid n$, \textit{there exists} some partition $\lambda^n_k$ such that $\gcd(\lambda^n_k)=m$. It is important to note, however, that there will still exist partitions for the same $n$ such that $\lambda^n_k\vdash$. For example, for 15 we have $15=5+10$, where $\gcd(5,10)=5$, but we also have $15=2+13$, for which $\gcd(2,13)=1$. (Note that $5\mid 15$ and $1\mid 15$.)
\end{remark}

\begin{lemma}\label{l:primes}
    $\lambda^p_k\vdash$ for any prime $p\geq k\geq2$.
\end{lemma}
\begin{proof}
    Let $p=\lambda_1+\cdots+\lambda_k$ and let $\gcd(\lambda_1,\dots,\lambda_k)=m$. There exist some positive integers $a_1,\dots,a_k$ such that $\lambda_i=a_im$ for all $1\leq i\leq k$. Then,
    \begin{align*}
        p&=\lambda_1+\cdots+\lambda_k\\
        &=a_1m+\cdots+a_km\\
        &=m(a_1+\cdots+a_k).
    \end{align*}
    We have $a_1+\cdots+a_k\geq2$ because $a_i\geq1$ for every $1\leq i\leq k$ and $k\geq2$, and therefore we need $m=1$ for $p$ to be prime, for we see that $\lambda^p_k\vdash$.
\end{proof}

\begin{corollary}\label{c:prime}
    Any integer $q$ is prime if and only if $\lambda_k^q\vdash$ for all $p\geq k\geq2$.
\end{corollary}

\begin{lemma}\label{l:all}
    For all $m\mid n$ and $1\leq k\leq n$ there exists $\gcd(\lambda^n_k)=m$. 
\end{lemma}
\begin{proof}
    If there exists some integer $a=a_1+\cdots+a_k$ such that $n=am$, then
    \begin{align*}
        am&=m(a_1+\cdots+a_k)\\
        n&=ma_1+\cdots+ma_k\\
        n&=\lambda_1+\cdots+\lambda_k,
    \end{align*}
    for we can clearly see that there exists $\gcd(\lambda_k^n)=m$ for some $m\mid n$. We use the restriction $1\leq k\leq n$ since $p(n,k)=0$ for all $k<1$ or $k>n$ (and thus we would have $\gcd(\lambda^n_k)=0\nmid n$).
    
    If there does not exist any integer $a$ which can be partitioned into $k$ parts such that $n=am$, then there does not exist $\gcd(\lambda^n_k)=m$. However, if $m\mid n$, then there exists some $k$ for which $\gcd(\lambda^n_k)=m$ since we have $a,m\geq1$.
\end{proof}

\begin{lemma}\label{l:iff}
    $\gcd(\lambda^n_k)=m$ for all $1\leq k\leq n$ if and only if $m\mid n$.
\end{lemma}
\begin{proof}
    Let $\gcd(\lambda^n_k)=m\nmid n$ and let $n=\lambda_1+\cdots+\lambda_k$. Since $m\nmid n$, there does not exist an integer $a=a_1+\cdots+a_k$ such that $n=am$ and 
    \begin{align*}
        am&=m(a_1+\cdots+a_k)\\
        n&=ma_1+\cdots+ma_k\\
        n&=\lambda_1+\cdots+\lambda_k,
    \end{align*}
    for we see that there are no partitions with $\gcd(\lambda^n_k)=m\nmid n$ and, thus, we only have partitions with $m\mid n$.
\end{proof}
Piecing everything together, we conclude by knowing that $\gcd(\lambda_k^n)=m$ if and only if $m\mid n$ for all $1\leq k\leq n$.
\end{proof}

\begin{remark}
    A way to rephrase \lref{iff} is: \textit{$\gcd(\lambda^n_k)\mid n$ for all $1\leq k\leq n$}.
\end{remark}

\begin{corollary}\label{c:general}
    $\lambda_k^n\vdash$ for all $m\mid n$ and $m\nmid\lambda_i\in\lambda_k^n$, that is, if at least one part of $\lambda_k^n$ does not share any common factor with $n$, then $\lambda_k^n\vdash$.
\end{corollary}

\begin{corollary}\label{c:formula}
    A partition of any nonnegative integer can only be relatively prime or not relatively prime. Furthermore, the number of partitions of a nonnegative integer is the number of relatively prime partitions plus the number of not-relatively prime partitions, that is, $p(n,k)=p_\Psi(n,k)+\Lambda(n,k)$ and $p(n)=p_\Psi(n)+\Lambda(n)$.
\end{corollary}

\section{Closed-form and recursive formula for the partition function}\label{s:formula}
Most methods of enumerating the number of (restricted or unrestricted) partitions of some nonnegative integer $n$ involve (1) analytic methods (especially through the theory of Modular Forms and manipulations of generating functions) or (2) visual combinatorial methods (especially Ferrers or Young diagrams). In this section, we present a method for representing partitions that differs from the usual method. We define partitions with $k$ parts, then organize them with their corresponding integer $n$, and then we define two formulas given a clear pattern that arises. We note that the formulas we seek may be inefficient, albeit combinatorially rich. 

Throughout, we will develop two formulas for $p(n,k)$: one in closed-form and a recursive version. The recurrence relation is easier to understand and prove by the method we will use, and then we can algebraically modify it to derive our closed-form formula. To give intuition for both expressions, we will explain them together and then formalize the results with Theorems \ref{th:recur1} and \ref{th:arithmetic_f}.

\begin{theorem}[Recursive formula]\label{th:recur1}
    In general, we have $p(n,1)=p(1)=1$ and for $n,k\geq2$ we have 
    \begin{align}\label{eq:rec_partk1}
        p(n,k)=\sum_{m=1}^{\left\lfloor\frac{n}{k}\right\rfloor}\sum_{\nu=1}^{k-1}p(n-km,\nu)
    \end{align}
    and
    \begin{align}\label{eq:rec_part1}
        p(n)=\sum_{k=1}^n\sum_{m=1}^{\left\lfloor\frac{n}{k}\right\rfloor}\sum_{\nu=1}^{k-1}p(n-km,\nu).
    \end{align}
\end{theorem}

\begin{theorem}[Closed-form formula]\label{th:arithmetic_f}
    We have $p(n,1)=1$, $p(n,2)=\left\lfloor\frac{n}{2}\right\rfloor$, for all $k\geq3$ we have
    \begin{align}\label{eq:sums_partk1}
        p(n,k)&=\sums_{\substack{t=0 \\ \left[k-t,\,1,\,0,\,\left\lfloor\frac{n-\sum_{j=k-t+1}^{k}jm_j}{k-t}\right\rfloor\right]}}^{k-3}\left\lfloor\frac{2+n-\sum_{j=3}^{k}jm_j}{2}\right\rfloor\\
        \notag&=\sum_{m_k=1}^{\left\lfloor\frac{n}{k}\right\rfloor}\sum_{m_{k-1}=0}^{\left\lfloor\frac{n-km_k}{k-1}\right\rfloor}\cdots\sum_{m_3=0}^{\left\lfloor\frac{n-4m_4-\cdots-(k-1)m_{k-1}-km_k}{3}\right\rfloor}\left\lfloor\frac{2+n-3m_3-\cdots-(k-1)m_{k-1}-km_k}{2}\right\rfloor
    \end{align}
    and for $n\geq0$ we have
    \begin{align}\label{eq:sums_part1}
        &p(n)=1+\left\lfloor\frac n2\right\rfloor+\sum_{k=3}^n\sums_{\substack{t=0 \\ \left[k-t,\,1,\,0,\,\left\lfloor\frac{n-\sum_{j=k-t+1}^{k}jm_j}{k-t}\right\rfloor\right]}}^{k-3}\left\lfloor\frac{2+n-\sum_{j=3}^{k}jm_j}{2}\right\rfloor\\
        \notag&=1+\left\lfloor\frac n2\right\rfloor+\sum_{m_3=1}^{\left\lfloor\frac{n}{3}\right\rfloor}\left\lfloor\frac{2+n-3m_3}{2}\right\rfloor+\cdots+\sum_{m_n=1}^{\left\lfloor\frac{n}{n}\right\rfloor}\cdots\sum_{m_3=0}^{\left\lfloor\frac{n-4m_4-\cdots-nm_n}{3}\right\rfloor}\left\lfloor\frac{2+n-3m_3-\cdots-nm_n}{2}\right\rfloor.
    \end{align}
\end{theorem}

\begin{remark}
    Due to \lref{primesp}, $p(n,1)=1$ for all $n$. This is the simplest formula for $p(n,k)$ and, as will be seen, $p(n,k)$ has more complicated formulas as the value of $k$ increases. 
\end{remark}

\subsection{Two parts}
In order to have control over the organization of partitions, we will use the notation of \tref{base2}, where we use partitions of the form $\lambda_2^n=\{\varSigma_1+c_1,~\varSigma_2+c_2\}$ and, for clarity, we do not include the brackets of the partitions. We denote this notation as \textit{base representation of $k=2$}. 
\begin{table}[t]\centering\begin{tabular}{c|c|c|cc}
    \multirow{2}{*}{$n$} & \multicolumn{4}{c}{$\lambda_2^n$} \\\cline{2-5}
    & $\cc0$ & $\cc1$ & $\cc2$ & $\cdots$ \\\hline
    $\varSigma_1+\varSigma_2$ & $\varSigma_1,~\varSigma_2$ && \\
    $\varSigma_1+\varSigma_2+1$ & $\varSigma_1,~\varSigma_2+1$ && \\
    $\varSigma_1+\varSigma_2+2$ & $\varSigma_1,~\varSigma_2+2$ & $\varSigma_1+1,~\varSigma_2+1$ & \\
    $\varSigma_1+\varSigma_2+3$ & $\varSigma_1,~\varSigma_2+3$ & $\varSigma_1+1,~\varSigma_2+2$ & \\
    $\varSigma_1+\varSigma_2+4$ & $\varSigma_1,~\varSigma_2+4$ & $\varSigma_1+1,~\varSigma_2+3$ & $\varSigma_1+2,~\varSigma_2+2$\\
    $\varSigma_1+\varSigma_2+5$ & $\varSigma_1,~\varSigma_2+5$ & $\varSigma_1+1,~\varSigma_2+4$ & $\varSigma_1+2,~\varSigma_2+3$\\
    \vdots&\vdots&\vdots&\vdots&$\ddots$
\end{tabular}\caption{Base representation of $k=2$}\label{t:base2}\end{table}
Note that we have placed a `$\cc{c_1}$' at the top of every group of columns that share this value of $c_1$. This notation will be useful for understanding our formula. We denote each group of columns that share this value of $c_1$ as $\cc{c_1}$-group, and we denote the number of partitions in each $\cc{c_1}$-group by $\#\cc{c_1}$. 

Let $\varSigma_1,\varSigma_2=1$, so that we get \tref{app2}. We denote this notation as \textit{applied representation of $k=2$}. 
\begin{table}[t]\centering\begin{tabular}{c|c|c|cc}
    \multirow{2}{*}{$n$} & \multicolumn{4}{c}{$\lambda_2^n$} \\\cline{2-5}
    & $\cc0$ & $\cc1$ & $\cc2$ & $\cdots$ \\\hline
    2 & 1,~1 && \\
    3 & 1,~2 && \\
    4 & 1,~3 & 2,~2 & \\
    5 & 1,~4 & 2,~3 & \\
    6 & 1,~5 & 2,~4 & 3,~3 \\
    7 & 1,~6 & 2,~5 & 3,~4 \\
    \vdots&\vdots&\vdots&\vdots&$\ddots$
\end{tabular}\caption{Applied representation of $k=2$}\label{t:app2}\end{table}
With it, it is easy to see that every partition of $n$ with $2$ parts can be represented with this format, given enough space. If we find a formula for the number of columns in any given layer $n$, we will consequently find a formula for $p(n,2)$. The sequence of values of $p(n,2)$, according to \tref{app2}, is $1,1,2,2,3,3,\dots$ for $n=2,3,4,\dots$ $\left\lfloor\frac n2\right\rfloor$ elegantly yields a result for this sequence, thus $p(n,2)=\left\lfloor\frac n2\right\rfloor$. (This result was found using other methods in, for example, \cite{hardy,honsberger}.)

For the recursive function, we may see the base and applied representations of $k=1$ (\tref{app1}) and understand that since we have $\left\lfloor\frac{n}{2}\right\rfloor$ columns in the base and applied representations of $k=2$ and since $p(n,1)=1$ for all $n\geq0$, we also have $p(n-2,1)$ repeated $\left\lfloor\frac{n}{2}\right\rfloor$ times.

\begin{table}[t]\centering
    \begin{subtable}[h]{0.5\textwidth}\centering
        \begin{tabular}{c|c}
            $n$ & $\lambda_1^n$ \\\hline
            $\varSigma_1$ & $\varSigma_1$ \\
            $\varSigma_1+1$ & $\varSigma_1+1$ \\
            $\varSigma_1+2$ & $\varSigma_1+2$ \\
            $\varSigma_1+3$ & $\varSigma_1+3$ \\
            $\varSigma_1+4$ & $\varSigma_1+4$ \\
            $\varSigma_1+5$ & $\varSigma_1+5$ \\
            $\varSigma_1+6$ & $\varSigma_1+6$ \\
            \vdots&\vdots
        \end{tabular}
        \caption{Base representation of $k=1$}
    \end{subtable}\hfill
    \begin{subtable}[h]{0.5\textwidth}\centering
        \begin{tabular}{c|c}
            $n$ & $\lambda_1^n$ \\\hline
            1 & 1 \\
            2 & 2 \\
            3 & 3 \\
            4 & 4 \\
            5 & 5 \\
            6 & 6 \\
            7 & 7 \\
            \vdots&\vdots
        \end{tabular}
        \caption{Applied representation of $k=1$}
    \end{subtable}
    \caption{Base and applied representations of $k=1$}\label{t:app1}
\end{table}

\subsection{Three parts}
For three parts, we have the base representation in \tref{base3}.
\begin{table}[t]\centering\small\begin{tabular}{c|c|c|c|c|c|c}
    \multirow{3}{*}{$n$} & \multicolumn{6}{c}{$\lambda_3^n$} \\\cline{2-7}
    & \multicolumn{2}{c|}{$\cc0$} & \multicolumn{2}{c|}{$\cc1$} & \multicolumn{1}{c}{$\cc2$} & $\cdots$ \\\cline{2-7}
    & {$\ccc0$} & $\cdots$ & {$\ccc1$} & $\cdots$ & {$\ccc2$} & $\cdots$ \\\hline
    $\varSigma_1+\varSigma_2+\varSigma_3$ & $\varSigma_1,~\varSigma_2,~\varSigma_3$ &&&& \\
    $\varSigma_1+\varSigma_2+\varSigma_3+1$ & $\varSigma_1,~\varSigma_2,~\varSigma_3+1$ &&&& \\
    $\varSigma_1+\varSigma_2+\varSigma_3+2$ & $\varSigma_1,~\varSigma_2,~\varSigma_3+2$ & $\cdots$ &&& \\
    $\varSigma_1+\varSigma_2+\varSigma_3+3$ & $\varSigma_1,~\varSigma_2,~\varSigma_3+3$ &  $\cdots$ & $\varSigma_1+1,~\varSigma_2+1,~\varSigma_3+1$ && \\
    $\varSigma_1+\varSigma_2+\varSigma_3+4$ & $\varSigma_1,~\varSigma_2,~\varSigma_3+4$ & $\cdots$ & $\varSigma_1+1,~\varSigma_2+1,~\varSigma_3+2$ && \\
    $\varSigma_1+\varSigma_2+\varSigma_3+5$ & $\varSigma_1,~\varSigma_2,~\varSigma_3+5$ & $\cdots$ & $\varSigma_1+1,~\varSigma_2+1,~\varSigma_3+3$ & $\cdots$ & \\
    $\varSigma_1+\varSigma_2+\varSigma_3+6$ & $\varSigma_1,~\varSigma_2,~\varSigma_3+6$ & $\cdots$ & $\varSigma_1+1,~\varSigma_2+1,~\varSigma_3+4$ & $\cdots$ & $\varSigma_1+2,~\varSigma_2+2,~\varSigma_3+2$ \\
    $\varSigma_1+\varSigma_2+\varSigma_3+7$ & $\varSigma_1,~\varSigma_2,~\varSigma_3+7$ & $\cdots$ & $\varSigma_1+1,~\varSigma_2+1,~\varSigma_3+5$ & $\cdots$ & $\varSigma_1+2,~\varSigma_2+2,~\varSigma_3+3$ \\
    $\varSigma_1+\varSigma_2+\varSigma_3+8$ & $\varSigma_1,~\varSigma_2,~\varSigma_3+8$ & $\cdots$ & $\varSigma_1+1,~\varSigma_2+1,~\varSigma_3+6$ & $\cdots$ & $\varSigma_1+2,~\varSigma_2+2,~\varSigma_3+4$ & $\cdots$ \\
    $\vdots$&$\vdots$&$\ddots$&$\vdots$&$\ddots$&$\vdots$&$\ddots$
\end{tabular}\caption{Base representation of $k=3$}\label{t:base3}\end{table}

\begin{remark}
    As before, we have partitions of the form $\lambda_3^n=\{\varSigma_1+c_1,~\varSigma_2+c_2,~\varSigma_3+c_3\}$. We note that, in general, for any $k\geq1$, we have partitions of the form $\lambda_k^n=\{\varSigma_1+c_1,\dots,~\varSigma_k+c_k\}$.
\end{remark}

Note that, as before, we have placed a `$\cc{c_1}$' at the top of every group of columns that share the same value of $c_1$, and now we have additionally placed a `$\ccc{c_2}$' at the top of every group of columns that share this value of $c_2$. We denote every $\ccc{c_2}$-group under $\cc{c_1}$ as $\ccc{c_2}_{\cc{c_1}}$-group, and the number of partitions in $\ccc{c_2}_{\cc{c_1}}$-group is $\#\ccc{c_2}_{\cc{c_1}}$. Since we use the convention $\lambda_1\leq\lambda_2\leq\cdots\leq\lambda_k$, then it follows that $c_1\leq c_2\leq\cdots\leq c_k$. (In this case, we have $c_1\leq c_2$.) This notation is clearer with the applied representation in \tref{app3}.
\begin{table}[t]\centering\small\begin{tabular}{c|c|c|c|c|c|c|c|c|c|c|c|c|c}
    \multirow{3}{*}{$n$} & \multicolumn{13}{c}{$\lambda_3^n$} \\\cline{2-14}
    & \multicolumn{6}{c|}{$\cc0$} & \multicolumn{4}{c|}{$\cc1$} & \multicolumn{2}{c}{$\cc2$} & $\cdots$ \\\cline{2-14}
    &$\ccc0$ & $\ccc1$ & $\ccc2$ & $\ccc3$ & $\ccc4$ & $\cdots$ & $\ccc1$ & $\ccc2$ & $\ccc3$ & $\cdots$ & $\ccc2$ & $\ccc3$ & $\cdots$ \\\hline
    3 & 1,~1,~1 &&&&&&&&&&& \\
    4 & 1,~1,~2 &&&&&&&&&&& \\
    5 & 1,~1,~3 & 1,~2,~2 &&&&&&&&&& \\
    6 & 1,~1,~4 & 1,~2,~3 &&&&& 2,~2,~2 &&&&& \\
    7 & 1,~1,~5 & 1,~2,~4 & 1,~3,~3 &&&& 2,~2,~3 &&&&& \\
    8 & 1,~1,~6 & 1,~2,~5 & 1,~3,~4 &&&& 2,~2,~4 & 2,~3,~3 &&&& \\
    9 & 1,~1,~7 & 1,~2,~6 & 1,~3,~5 & 1,~4,~4 &&& 2,~2,~5 & 2,~3,~4 &&& 3,~3,~3 & \\
    10 & 1,~1,~8 & 1,~2,~7 & 1,~3,~6 & 1,~4,~5 &&& 2,~2,~6 & 2,~3,~5 & 2,~4,~4 && 3,~3,~4 & \\
    11 & 1,~1,~9 & 1,~2,~8 & 1,~3,~7 & 1,~4,~6 & 1,~5,~5 && 2,~2,~7 & 2,~3,~6 & 2,~4,~5 && 3,~3,~5 & 3,~4,~4 \\
    $\vdots$&$\vdots$&$\vdots$&$\vdots$&$\vdots$&$\vdots$&$\ddots$&$\vdots$&$\vdots$&$\vdots$&$\ddots$&$\vdots$&$\vdots$&$\ddots$
\end{tabular}\caption{Applied representation of $k=3$}\label{t:app3}\end{table}

Regardless of the actual values of $n$ and those of the parts, we can observe that each $\cc{c_1}$-group has the same pattern that we had in \tref{app2}. Note that every partition in $\cc0$-group, for example, is `shifted by $+1$' compared to $p(n,2)$ since we consider an extra part plus the two parts that we already had in \tref{app2}. We can thus formulate that
\begin{align*}
    \#\cc0=\left\lfloor\frac{n-1}{2}\right\rfloor=\left\lfloor\frac{2+n-3}{2}\right\rfloor.
\end{align*}

For $\cc1$-group, note that every partition is `shifted by $+3$' when compared with $\cc0$-group since the parts of the first partition of $\cc1$-group ($\{2,2,2\}$) are the parts of the first partition of $\cc0$-group ($\{1,1,1\}$) plus 1 and the rest of the partitions follow in the same manner. Thus, we can formulate that
\begin{align*}
    \#\cc1=\left\lfloor\frac{2+n-3-3}{2}\right\rfloor=\left\lfloor\frac{2+n-6}{2}\right\rfloor.
\end{align*}

For $\cc2$-group the same occurs: all the parts of the first partition of $\cc2$-group ($\{3,3,3\}$) are equivalent to the parts of the first partition of $\cc1$-group ($\{2,2,2\}$) plus 1 and the rest of the partitions follow in the same manner. Thus, we can formulate that 
\begin{align*}
    \#\cc2=\left\lfloor\frac{2+n-6-3}{2}\right\rfloor=\left\lfloor\frac{2+n-9}{2}\right\rfloor.
\end{align*}

Note that this process will everlastingly occur for every $\cc{c_1}$-group, for we can generalize
\begin{align*}
    \#\cc{c_1}&=\left\lfloor\frac{2+n-3(c_1+1)}{2}\right\rfloor.
\end{align*}
The number of partitions of $n$ with 3 parts is the sum of all $\#\cc{c_1}$ for the corresponding value of $n$, that is,
\begin{align}
    \notag p(n,3)&=\sum_{c_1=0}^{\ell-1}\#\cc{c_1}\\
    \notag&=\sum_{c_1=0}^{\ell-1}\left\lfloor\frac{2+n-3(c_1+1)}{2}\right\rfloor\\
    \label{eq:pn3}&=\sum_{m=1}^\ell\left\lfloor\frac{2+n-3m}{2}\right\rfloor.
\end{align}
We require an upper bound $\ell$ to avoid getting negative values. Since every $n$ has partitions in exactly $\left\lfloor\frac{n}{3}\right\rfloor$ columns, $\ell=\left\lfloor\frac{n}{3}\right\rfloor$. (We can also set the summand of \eqref{eq:pn3} to be $\geq1$ and solve the inequality to find the same result.)

The recursive formula for $k=3$ can be explained as follows. Since the first two parts of the partitions in $\ccc0_{\cc0}$-group are equal to 1, when we subtract $n-3$, we are essentially subtracting 1 from every part, and thus we are left with only one remaining part, for we see that $\#\ccc0_{\cc0}=p(n-3,1)$. Moreover, since the rest of the partitions in $\cc0$-group have $\lambda_1=1$ (and $\lambda_2,\lambda_3>1$), when we repeat this process, subtracting $n-3$, we get $p(n-3,2)$ partitions (since $\lambda_1$ is not considered anymore and $\lambda_2,\lambda_3\geq1$), for, in general,
\begin{align*}
    \#\cc0=p(n-3,1)+p(n-3,2).
\end{align*}

When we subtract 3 from the partitions in $\cc1$-group, we essentially subtract 1 from every part of these partitions, for we see that we arrive again at what we had with $\cc0$-group. From here, it is clear that $\#\cc1=p(n-3-3,1)+p(n-3-3,2)=p(n-6,1)+p(n-6,2)$. Moving onto $\cc2$-group, the same recursion occurs when we subtract 3 from every partition, and it is easy to show that (see below for a rigorous proof), in general, 
\begin{align*}
    \#\cc{c_1}&=p(n-3(c_1+1),1)+p(n-3(c_1+1),2)
\end{align*}
and, moreover, 
\begin{align*}
    p(n,3)&=\sum_{c_1=0}^{\ell-1}\#\cc{c_1}\\
    &=\sum_{c_1=0}^{\ell-1}p(n-3(c_1+1),1)+p(n-3(c_1+1),2)\\
    &=\sum_{m=1}^{\ell}p(n-3m,1)+p(n-3m,2).
\end{align*}

Lastly, since every $n$ has partitions in exactly $\left\lfloor\frac{n}{3}\right\rfloor$ columns, consequently $\ell=\left\lfloor\frac{n}{3}\right\rfloor$.

\subsection{Four parts}
\begin{table}[t]\centering\scriptsize\setlength{\tabcolsep}{0.21em}\begin{tabular}{c|c|c|c|c|c|c}
    \multirow{4}{*}{$n$} & \multicolumn{6}{c}{$\lambda_4^n$} \\\cline{2-7}
    & \multicolumn{2}{c|}{$\cc{0}$} & \multicolumn{2}{c|}{$\cc{1}$} & \multicolumn{1}{c}{$\cc{2}$} & $\cdots$ \\\cline{2-7}
    & {$\ccc0$} & $\cdots$ & {$\ccc1$} & $\cdots$ & {$\ccc2$} & $\cdots$ \\\cline{2-7}
    & {$\cccc0$} & {$\cdots$} & {$\cccc1$} & {$\cdots$} & {$\cccc2$} & {$\cdots$} \\\hline
    $\varSigma_1+\varSigma_2+\varSigma_3+\varSigma_4$ & $\varSigma_1,~\varSigma_2,~\varSigma_3,~\varSigma_4$ &&&& \\
    $\varSigma_1+\varSigma_2+\varSigma_3+\varSigma_4+1$ & $\varSigma_1,~\varSigma_2,~\varSigma_3,~\varSigma_4+1$ &&&& \\
    $\varSigma_1+\varSigma_2+\varSigma_3+\varSigma_4+2$ & $\varSigma_1,~\varSigma_2,~\varSigma_3,~\varSigma_4+2$ & $\cdots$ &&& \\
    $\varSigma_1+\varSigma_2+\varSigma_3+\varSigma_4+3$ & $\varSigma_1,~\varSigma_2,~\varSigma_3,~\varSigma_4+3$ & $\cdots$ &&& \\
    $\varSigma_1+\varSigma_2+\varSigma_3+\varSigma_4+4$ & $\varSigma_1,~\varSigma_2,~\varSigma_3,~\varSigma_4+4$ & $\cdots$ & $\varSigma_1+1,~\varSigma_2+1,~\varSigma_3+1,~\varSigma_4+1$ && \\
    $\varSigma_1+\varSigma_2+\varSigma_3+\varSigma_4+5$ & $\varSigma_1,~\varSigma_2,~\varSigma_3,~\varSigma_4+5$ & $\cdots$ & $\varSigma_1+1,~\varSigma_2+1,~\varSigma_3+1,~\varSigma_4+2$ && \\
    $\varSigma_1+\varSigma_2+\varSigma_3+\varSigma_4+6$ & $\varSigma_1,~\varSigma_2,~\varSigma_3,~\varSigma_4+6$ & $\cdots$ & $\varSigma_1+1,~\varSigma_2+1,~\varSigma_3+1,~\varSigma_4+3$ & $\cdots$ & \\
    $\varSigma_1+\varSigma_2+\varSigma_3+\varSigma_4+7$ & $\varSigma_1,~\varSigma_2,~\varSigma_3,~\varSigma_4+7$ & $\cdots$ & $\varSigma_1+1,~\varSigma_2+1,~\varSigma_3+1,~\varSigma_4+4$ & $\cdots$ && \\
    $\varSigma_1+\varSigma_2+\varSigma_3+\varSigma_4+8$ & $\varSigma_1,~\varSigma_2,~\varSigma_3,~\varSigma_4+8$ & $\cdots$ & $\varSigma_1+1,~\varSigma_2+1,~\varSigma_3+1,~\varSigma_4+5$ & $\cdots$ &  $\varSigma_1+2,~\varSigma_2+2,~\varSigma_3+2,~\varSigma_4+2$ & \\
    $\vdots$&$\vdots$&$\ddots$&$\vdots$&$\ddots$&$\vdots$&$\ddots$
\end{tabular}\caption{Base representation of $k=4$}\label{t:base4}\end{table}

\begin{table}[t]\centering\scriptsize\setlength{\tabcolsep}{0.59em}\begin{tabular}{c|c|c|c|c|c|c|c|c|c|c|c|c|c}
    \multirow{4}{*}{$n$} & \multicolumn{13}{c}{$\lambda_4^n$} \\\cline{2-14}
    & \multicolumn{6}{c|}{$\cc{0}$} & \multicolumn{5}{c|}{$\cc{1}$} & \multicolumn{1}{c}{$\cc{2}$} & $\cdots$ \\\cline{2-14}
    & \multicolumn{3}{c|}{$\ccc0$} & \multicolumn{2}{c|}{$\ccc1$} & $\cdots$ & \multicolumn{3}{c|}{$\ccc1$} & {$\ccc2$} & $\cdots$ & {$\ccc2$} & $\cdots$ \\\cline{2-14}
    & {$\cccc0$} & {$\cccc1$} & {$\cdots$} & {$\cccc1$} & {$\cccc2$} & {$\cdots$} & {$\cccc1$} & {$\cccc2$} & $\cdots$ & {$\cccc2$} & $\cdots$ & {$\cccc2$} & {$\cdots$} \\\hline
    4 & 1,~1,~1,~1 &&&&&&&&&&&& \\
    5 & 1,~1,~1,~2 &&&&&&&&&&&& \\
    6 & 1,~1,~1,~3 & 1,~1,~2,~2 &&&&&&&&&&& \\
    7 & 1,~1,~1,~4 & 1,~1,~2,~3 && 1,~2,~2,~2 &&&&&&&&& \\
    8 & 1,~1,~1,~5 & 1,~1,~2,~4 & $\cdots$ & 1,~2,~2,~3 &&& 2,~2,~2,~2 &&&&&& \\
    9 & 1,~1,~1,~6 & 1,~1,~2,~5 & $\cdots$ & 1,~2,~2,~4 & 1,~2,~3,~3 && 2,~2,~2,~3 &&&&&& \\
    10 & 1,~1,~1,~7 & 1,~1,~2,~6 & $\cdots$ & 1,~2,~2,~5 & 1,~2,~3,~4 && 2,~2,~2,~4 & 2,~2,~3,~3 &&&& \\
    11 & 1,~1,~1,~8 & 1,~1,~2,~7 & $\cdots$ & 1,~2,~2,~6 & 1,~2,~3,~5 & $\cdots$ & 2,~2,~2,~5 & 2,~2,~3,~4 && 2,~3,~3,~3 && \\
    12 & 1,~1,~1,~9 & 1,~1,~2,~8 & $\cdots$ & 1,~2,~2,~7 & 1,~2,~3,~6 & $\cdots$ & 2,~2,~2,~6 & 2,~2,~3,~5 & $\cdots$ & 2,~3,~3,~4 && 3,~3,~3,~3 \\
    $\vdots$&$\vdots$&$\vdots$&$\ddots$&$\vdots$&$\vdots$&$\ddots$&$\vdots$&$\vdots$&$\ddots$&$\vdots$&$\ddots$&$\vdots$&$\ddots$
\end{tabular}\caption{Applied representation of $k=4$}\label{t:app4}\end{table}

For four parts, we have the base representation in \tref{base4} and the applied representation in \tref{app4}. We see that $\ccc0_{\cc0}$-group is `shifted by $+1$' compared to $\cc0$-group of \tref{app3} since we consider an additional part plus the three parts that we already had in that group, for we can generalize
\begin{align*}
    \#\ccc0_{\cc0}=\left\lfloor\frac{2+n-3-1}{2}\right\rfloor=\left\lfloor\frac{2+n-4}{2}\right\rfloor.
\end{align*}

Similar to how we worked with $k=3$ parts, note that every partition in $\ccc1_{\cc0}$-group is `shifted by $+3$' when compared with $\ccc0_{\cc0}$-group since the last three parts of the first partition of $\ccc1_{\cc0}$-group ($\{1,2,2,2\}$) are the same last three parts of the first partition of $\ccc0_{\cc0}$-group ($\{1,1,1,1\}$) plus 1, the first part is the same, and the rest of the partitions follow in the same manner, for we can formulate that
\begin{align*}
    \#\ccc1_{\cc0}=\left\lfloor\frac{2+n-3-4}{2}\right\rfloor.
\end{align*}

For $\ccc2_{\cc0}$-group a similar thing occurs: every partition is `shifted by $+3$' compared to $\ccc1_{\cc0}$-group since the last three parts of the first partition of $\ccc2_{\cc0}$-group ($\{1,3,3,3\}$) are equivalent to the last three parts of the first partition of $\ccc1_{\cc0}$-group ($\{1,2,2,2\}$) plus 1, the first part is the same, and the rest of the partitions follow in the same manner, for we can formulate that 
\begin{align*}
    \#\ccc2_{\cc0}=\left\lfloor\frac{2+n-3-3-4}{2}\right\rfloor=\left\lfloor\frac{2+n-6-4}{2}\right\rfloor.
\end{align*}

Note the generalization
\begin{align*}
    \#\ccc{c_2}_{\cc0}&=\left\lfloor\frac{2+n-3c_2-4}{2}\right\rfloor.
\end{align*}
We can generalize these results to
\begin{align*}
    \#\cc0&=\sum_{c_2=0}^\ell\#\ccc{c_2}_{\cc{0}}\\
    &=\sum_{c_2=0}^\ell\left\lfloor\frac{2+n-3c_2-4}{2}\right\rfloor\\
    &=\sum_{m=0}^\ell\left\lfloor\frac{2+n-3m-4}{2}\right\rfloor.
\end{align*}
Since summing terms which are equal to 0 does not affect the result, we have
\begin{align*}
    \left\lfloor\frac{2+n-3m-4}2\right\rfloor&\geq1\\
    \frac{2+n-3m-4}2&\geq1\\
    n-3m-4&\geq0\\
    m&\leq\frac{n-4}3\\
    \ell&=\left\lfloor\frac{n-4}{3}\right\rfloor.
\end{align*}

Now note that $\ccc1_{\cc1}$-group is `shifted by 4' compared to $\ccc0_{\cc0}$-group since all the parts of the first partition of $\ccc1_{\cc1}$-group ($\{2,2,2,2\}$) are equivalent to all the parts of the first partition of $\ccc0_{\cc0}$-group ($\{1,1,1,1\}$) plus 1 and the rest of the partitions follow in the same manner, for we can formulate
\begin{align*}
    \#\ccc1_{\cc1}&=\left\lfloor\frac{2+n-4-4}{2}\right\rfloor=\left\lfloor\frac{2+n-8}{2}\right\rfloor.
\end{align*}

Moreover, similar to before, $\ccc2_{\cc1}$-group is `shifted by 3' compared to $\ccc1_{\cc1}$-group since the last three parts of the first partition of $\ccc2_{\cc1}$-group ($\{2,3,3,3\}$) are the same as the last three parts of the first partition of $\ccc1_{\cc1}$-group ($\{2,2,2,2\}$) plus 1, the first part is the same, and the rest of the partitions follow in the same manner, for we can formulate that
\begin{align*}
    \#\ccc2_{\cc1}=\left\lfloor\frac{2+n-3-8}{2}\right\rfloor.
\end{align*}

This pattern will continue indefinitely by
\begin{align*}
    \#\ccc{c_2}_{\cc1}=\left\lfloor\frac{2+n-3(c_2-1)-8}{2}\right\rfloor,
\end{align*}
for we can see that
\begin{align*}
    \#\cc1&=\sum_{c_2=1}^{\ell+1}\#\ccc{c_2}_{\cc1}\\
    &=\sum_{c_2=1}^{\ell+1}\left\lfloor\frac{2+n-3(c_2-1)-8}{2}\right\rfloor\\
    &=\sum_{m=0}^\ell\left\lfloor\frac{2+n-3m-8}{2}\right\rfloor
\end{align*}
and
\begin{align*}
    \left\lfloor\frac{2+n-3m-8}2\right\rfloor&\geq1\\
    \frac{2+n-3m-8}2&\geq1\\
    n-3m-8&\geq0\\
    m&\leq\frac{n-8}3\\
    \ell&=\left\lfloor\frac{n-8}{3}\right\rfloor.
\end{align*}

For $\cc2$-group the same pattern arises, for we can see that (continuing the same process)
\begin{align*}
    \#\cc2&=\sum_{c_2=2}^{\ell+2}\left\lfloor\frac{2+n-3(c_2-2)-12}{2}\right\rfloor\\
    &=\sum_{m=0}^\ell\left\lfloor\frac{2+n-3m-12}{2}\right\rfloor
\end{align*}
and
\begin{align*}
    \left\lfloor\frac{2+n-3m-12}2\right\rfloor&\geq1\\
    \frac{2+n-3m-12}2&\geq1\\
    n-3m-12&\geq0\\
    m&\leq\frac{n-12}3\\
    \ell&=\left\lfloor\frac{n-12}{3}\right\rfloor.
\end{align*}

We can therefore note the pattern
\begin{align*}
    p(n,4)&=\#\cc0+\#\cc1+\#\cc2+\cdots\\
    &=\sum_{m=0}^{\left\lfloor\frac{n-4}{3}\right\rfloor}\left\lfloor\frac{2+n-3m-4}{2}\right\rfloor+\sum_{m=0}^{\left\lfloor\frac{n-8}{3}\right\rfloor}\left\lfloor\frac{2+n-3m-8}{2}\right\rfloor+\sum_{m=0}^{\left\lfloor\frac{n-12}{3}\right\rfloor}\left\lfloor\frac{2+n-3m-12}{2}\right\rfloor+\cdots\\
    &=\sum_{m_4=1}^\ell\sum_{m_3=0}^{\left\lfloor\frac{n-4m_4}{3}\right\rfloor}\left\lfloor\frac{2+n-3m_3-4m_4}{2}\right\rfloor.
\end{align*}
Since every $n$ has partitions in exactly $\left\lfloor\frac{n}{4}\right\rfloor$ columns, consequently $\ell=\left\lfloor\frac{n}{4}\right\rfloor$.

Similar to before, the recursive formula for $k=4$ can be explained as follows. Since the first three parts of the partitions of $\cccc0_{\ccc0_{\cc{0}}}$-group are equal to 1, when we subtract $n-4$ we see that we essentially subtract 1 from each part, and thus we see that we have $\lambda_1=\lambda_2=\lambda_3=0$ and $\lambda_4\geq0$, for we have $\#\cccc0_{\ccc0_{\cc{0}}}=p(n-4,1)$. Thereafter, the first two parts in $\cccc{c_3}_{\ccc0_{\cc{0}}}$-group ($c_3>0$) are equal to 1, for we can subtract $n-4$ to see that $\sum_{c_3>0}\#\cccc{c_3}_{\ccc0_{\cc{0}}}=p(n-4,2)$. Lastly, repeating the same procedure, we may find that $\sum_{c_2>0}\#\ccc{c_2}_{\cc0}=p(n-4,3)$, for we can generalize
\begin{align*}
    \#{\cc{0}}=p(n-4,1)+p(n-4,2)+p(n-4,3).
\end{align*}

If we subtract 4 from all partitions in $\cc1$-group, it is easy to see that we get the number of partitions in $\cc0$-group, that is,
\begin{align*}
    \#\cc1&=p(n-4-4,1)+p(n-4-4,2)+p(n-4-4,3)\\
    &=p(n-8,1)+p(n-8,2)+p(n-8,3).
\end{align*}
We can see that this will occur for all values of $c_1\geq0$, for we can generalize
\begin{align*}
    \#\cc{c_1}&=p(n-4(c_1+1),1)+p(n-4(c_1+1),2)+p(n-4(c_1+1),3)\\
    \implies p(n,4)&=\sum_{c_1=0}^{\ell-1}\#\cc{c_1}\\
    &=\sum_{c_1=0}^{\ell-1}p(n-4(c_1+1),1)+p(n-4(c_1+1),2)+p(n-4(c_1+1),3)\\
    &=\sum_{m=1}^{\ell}p(n-4m,1)+p(n-4m,2)+p(n-4m,3).
\end{align*}
Since every $n$ has partitions in exactly $\left\lfloor\frac{n}{4}\right\rfloor$ columns, consequently $\ell=\left\lfloor\frac{n}{4}\right\rfloor$.

\subsection{Generalization}
Let $B_k$ be the base representation of some arbitrary value of $k\geq1$. Let $\cc{c_i}^i$ be the generalization of the notation we have used throughout this section to denote the different groups of partitions in our base and applied representations, defined recursively as
\begin{align*}
    \cc{c_1}^1=\cc{c_1};\quad\cc{c_i}^{i}=\cc{\cc{c_i}^{i-1}}\quad(i>1).
\end{align*}
Let $\cc{c_i}^i_{\cc{c_{i-1}}^{i-1}_{\cdots}}$ be equal to the notation $\cc{c_1:c_2:\cdots:c_{i-2}:c_{i-1}:c_i}$ (with $i<k$), and we will write $\cc{c_1:\cdots:c_i}_k$-group to specify that $\cc{c_1:\cdots:c_i}$-group belongs to $B_k$. With this, we may write, for example, $\#\cc{0:2:4}_{11}$ as a simplification of $\#\cccc4_{\ccc2_{\cc0}}$ of $B_{11}$, and $\cc6^5$ as a simplification of $\cccc{\ccc6}$.

\begin{remark}
    Before proceeding with the generalization of $p(n,k)$, it is important to show \thref{representations}.
\end{remark}

\begin{theorem}\label{th:representations}
    All partitions with $k$ parts of some integer $n\geq k$ can be found in $B_k$.
\end{theorem}

\begin{proof}
    All partitions in $B_k$ have the form $\{\varSigma_1+c_1,~\varSigma_2+c_2,\dots,~\varSigma_k+c_k\}$. Fix $c_1,c_2,\dots,c_{k-1}$ to be some arbitrary constant. If
    \begin{align}\label{eq:B_k}
        \sum_{i=1}^k\varSigma_i+c_i>n,
    \end{align}
    then there will be zero partitions of $n$ in $|c_1:c_2:\cdots:c_{k-1}|$-group (regardless of the value of $k$). Otherwise, consider \[m=n-c_1-c_2-\cdots-c_{k-1}\] and $k-(k-1)=1$ part. It is clear that $B_1$ contains all partitions with 1 part of the integers $m\geq1$. 

    Now, fix $c_1,c_2,\dots,c_{k-j}$ to be some arbitrary constant. Consider the integer \[r=n-c_1-c_2-\cdots-c_{k-j}\] and $k-(k-j)=j$ parts. Repeat this process until we are left with only one part. Then, clearly, all partitions of $r$ with given $c_1,c_2,\dots,c_{k-1}$ values will have their corresponding partition in $B_1$ (as long as \eqref{eq:B_k} holds, otherwise there will not be partitions of $r$ in $B_1$). Change the values of $c_1,c_2,\dots,c_{k-1}$ sufficiently, and note that this implies that any partition with $k$ parts of some integer $n\geq k$ can be found in $B_k$.
\end{proof}

\begin{proof}[Proof of \thref{recur1}]
    To show \eqref{eq:rec_partk1} and \eqref{eq:rec_part1}, we shall show the following:
    \begin{enumerate}
        \item[(i)] If we let $c_1=c_2=\cdots=c_{k-j-1}=0$ and $1\leq c_{k-j}\leq c_{k-j+1}\leq\cdots\leq c_{k-1}\leq c_{k}$ for $0\leq j\leq k-2$, then $\#|0:\cdots:0:c_{k-j}:c_{k-j+1}:\cdots:c_{k-2}:c_{k-1}|_k=p(n-k,j+1)$. 
        \item[(ii)] $\#\cc{c_1}_k=\sum_{\nu=1}^{k-1}p(n-k(c_1+1),\nu)$.
    \end{enumerate}

    \begin{proof}[Proof of (i)]
        We have $\lambda_i=\varSigma_i+c_i=1+c_i$. Since $c_i=0$ for $1\leq i\leq k-j-1$, then $\lambda_i=1$. Moreover, we have
       \begin{align}\label{eq:xparts}
            n&=\underbrace{\lambda_1+\lambda_2+\cdots+\lambda_{k-j-1}}_{k-j-1\text{ parts}}+\underbrace{\lambda_{k-j}+\cdots+\lambda_{k-1}+\lambda_k}_{x\text{ parts}}\\
            \notag&=\underbrace{1+1+\cdots+1}_{k-j-1\text{ parts}}+\underbrace{\varSigma_{k-j}+c_{k-j}+\cdots+\varSigma_{k-1}+c_{k-1}+\varSigma_k+c_k}_{x\text{ parts}}\\
            \notag&=k-j-1+x+c_{k-j}+\cdots+c_{k-1}+c_k.
        \end{align}
        Since the right-hand side of \eqref{eq:xparts} has $k$ parts, then
        \begin{align*}
            k-j-1+x&=k\\
            x&=j+1.
        \end{align*}
        Thereafter, we have
        \begin{align}
            \notag n&=k-j-1+x+c_{k-j}+\cdots+c_{k-1}+c_k\\
            \notag&=k+c_{k-j}+\cdots+c_{k-1}+c_k\\
            \label{eq:j+1parts}n-k&=\underbrace{c_{k-j}+\cdots+c_{k-1}+c_k}_{j+1\text{ parts}}.
        \end{align}
        All the parts $\lambda_1,\lambda_2,\dots,\lambda_{k-j-1}$ are fixed to be equal to 1, and thus we can consider them as constants, for which $\#|0:\cdots:0:c_{k-j}:\cdots:c_{k-1}|_k$ only depends on the parts $\lambda_{k-j},\lambda_{k-j+1},\dots,\lambda_{k-1},\lambda_k$. If we subtract $n-k$, we `eliminate' the parts $\lambda_{1},\lambda_{2},\dots,\lambda_{k-j-1}$, and thus we are left with $\#|0:\cdots:0:c_{k-j}:\cdots:c_{k-1}|_k=p(n-k,j+1)$ (where we consider $j+1$ parts since this is the amount of parts on the right-hand side of \eqref{eq:j+1parts}).
    \end{proof}
    \begin{remark}
        In (i) we use the restriction $0\leq j\leq k-2$ since when $j=0$ we have $c_1=c_2=\cdots=c_{k-1}=0$ and $1\leq c_k$\footnote{\label{fn:p(0)}$k=1$ is the only value of $k$ for which it is possible that $c_k=c_1=0$ since, even though $n-k=0$ when $n=k$, $p(0,1)=1$ (unlike with other values of $k$).}, and when $j=k-2$, we have $c_1=0$ and $1\leq c_2\leq c_3\leq\cdots\leq c_k$. In any case, it is important that we have $c_1=0$ since we will deal with $c_1\geq0$ in (ii).
    \end{remark}

    \begin{proof}[Proof of (ii)]
        To set the grounds for all values of $c_1\geq0$, we must start with $c_1=0$. Due to \thref{representations} and given that $0\leq j\leq k-2$, we may formulate
        \begin{align}
            \label{eq:base1(ii)}\#\cc0_k&=\sum_{j=0}^{k-2}p(n-k,j+1)\\
            \label{eq:base(ii)}&=\sum_{\nu=1}^{k-1}p(n-k,\nu).
        \end{align}
        
        Moreover, let $\lambda_k^n=\{\lambda_1,\dots,\lambda_k\}$ and consider $\lambda_1>1$. We may see subtracting $n-k$ as subtracting 1 from every part in $\lambda_k^n$ since we have exactly $k$ parts and $k=\sum_{i=1}^k1$ implies
        \begin{align*}
            \sum_{i=1}^k\lambda_i-1=n-k.
        \end{align*}
        Thus, we must subtract $k$ from $n$ exactly $\lambda_1-1=c_1$ times to get $\lambda_1=1$ and be able to use \eqref{eq:base(ii)}. We thus have $\#\cc0_k=\sum_{\nu=1}^{k-1}p(n-k,\nu)$, $\#\cc1_k=\sum_{\nu=1}^{k-1}p(n-2k,\nu)$, $\#\cc2_k=\sum_{\nu=1}^{k-1}p(n-3k,\nu)$, etc., for we see that, in general,
        \begin{align}\label{eq:c1k}
            \#\cc{c_1}_k=\sum_{\nu=1}^{k-1}p(n-k(c_1+1),\nu).
        \end{align}
    \end{proof}
    
    Since $n$ has partitions in exactly $\left\lfloor\frac{n}{k}\right\rfloor$ $\cc{c_1}_k$-groups, we have
        \begin{align*}
            p(n,k)=\sum_{c_1=0}^{\left\lfloor\frac{n}{k}\right\rfloor-1}\#\cc{c_1}_k=\sum_{c_1=0}^{\left\lfloor\frac{n}{k}\right\rfloor-1}\sum_{\nu=1}^{k-1}p(n-k(c_1+1),\nu)=\sum_{m=1}^{\left\lfloor\frac{n}{k}\right\rfloor}\sum_{\nu=1}^{k-1}p(n-km,\nu)
        \end{align*}
        and
        \begin{align*}
            p(n)=\sum_{k=1}^np(n,k)=\sum_{k=1}^n\sum_{m=1}^{\left\lfloor\frac{n}{k}\right\rfloor}\sum_{\nu=1}^{k-1}p(n-km,\nu).
        \end{align*}
\end{proof}

\begin{proof}[Proof of \thref{arithmetic_f}]
Since we have shown that \eqref{eq:rec_partk1} is true, we can check that \eqref{eq:sums_partk1} is equal to \eqref{eq:rec_partk1}, and then manipulate it to prove \thref{arithmetic_f}. 
We have $p(n,1)=p(1)=1$, and we can inductively see that
\begin{align*}
    p(n,2)&=\left\lfloor\frac{n}{2}\right\rfloor=\sum_{m=1}^{\left\lfloor\frac{n}{2}\right\rfloor}1=\sum_{m=1}^{\left\lfloor\frac{n}{2}\right\rfloor}p(n-2m,1)=\sum_{m=1}^{\left\lfloor\frac{n}{2}\right\rfloor}\sum_{\nu=1}^{2-1}p(n-2m,\nu),\\
    p(n,3)&=\sum_{m=1}^{\left\lfloor\frac{n}{3}\right\rfloor}\left\lfloor\frac{2+n-3m}{2}\right\rfloor=\sum_{m=1}^{\left\lfloor\frac{n}{3}\right\rfloor}p(n-3m,1)+p(n-3m,2)=\sum_{m=1}^{\left\lfloor\frac{n}{3}\right\rfloor}\sum_{\nu=1}^{3-1}p(n-3m,\nu),
\end{align*}
and, for some $k\geq3$,
\begin{align*}
    &p(n,k)=\sums_{\substack{t=0 \\ \left[k-t,\,1,\,0,\,\left\lfloor\frac{n-\sum_{j=k-t+1}^{k}jm_j}{k-t}\right\rfloor\right]}}^{k-3}\left\lfloor\frac{2+n-\sum_{j=3}^{k}jm_j}{2}\right\rfloor\\
    &=\sum_{m_k=1}^{\left\lfloor\frac{n}{k}\right\rfloor}
    \sum_{m_{k-1}=0}^{\left\lfloor\frac{n-km_k}{k-1}\right\rfloor}
    \cdots\sum_{m_3=0}^{\left\lfloor\frac{n-4m_4-\cdots-(k-1)m_{k-1}-km_k}{3}\right\rfloor}
    \left\lfloor\frac{2+n-3m_3-\cdots-(k-1)m_{k-1}-km_k}{2}\right\rfloor\\
    &=\sum_{m_k=1}^{\left\lfloor\frac{n}{k}\right\rfloor}
    \sum_{m_{k-2}=0}^{\left\lfloor\frac{n-(k-1)(0)-km_k}{k-2}\right\rfloor}\cdots
    \sum_{m_3=0}^{\left\lfloor\frac{n-4m_4-\cdots-(k-1)(0)-km_k}{3}\right\rfloor}
    \left\lfloor\frac{2+n-3m_3-\cdots-(k-1)(0)-km_k}{2}\right\rfloor
    \\&\quad+\sum_{m_{k-1}=1}^{\left\lfloor\frac{n-km_k}{k-1}\right\rfloor}
    \cdots\sum_{m_3=0}^{\left\lfloor\frac{n-4m_4-\cdots-(k-1)m_{k-1}-km_k}{3}\right\rfloor}
    \left\lfloor\frac{2+n-3m_3-\cdots-(k-1)m_{k-1}-km_k}{2}\right\rfloor\\
    &=\sum_{m_k=1}^{\left\lfloor\frac{n}{k}\right\rfloor}
    \sum_{m_{k-3}=0}^{\left\lfloor\frac{n-(k-2)(0)-km_k}{k-3}\right\rfloor}\cdots
    \sum_{m_3=0}^{\left\lfloor\frac{n-4m_4-\cdots-(k-2)(0)-km_k}{3}\right\rfloor}
    \left\lfloor\frac{2+n-3m_3-\cdots-(k-2)(0)-km_k}{2}\right\rfloor
    \\&\quad+\sum_{m_{k-2}=1}^{\left\lfloor\frac{n-km_k}{k-2}\right\rfloor}
    \sum_{m_3=0}^{\left\lfloor\frac{n-4m_4-\cdots-(k-2)m_{k-2}-km_k}{3}\right\rfloor}
    \left\lfloor\frac{2+n-3m_3-\cdots-(k-2)m_{k-2}-km_k}{2}\right\rfloor
    +p(n-km_k,k-1)\\
    &=\cdots=\sum_{m_k=1}^{\left\lfloor\frac{n}{k}\right\rfloor}\left\lfloor\frac{2+n-3(0)-km_k}{2}\right\rfloor+
    \sum_{m_3=1}^{\left\lfloor\frac{n-km_k}{3}\right\rfloor}
    \left\lfloor\frac{2+n-3m_3-km_k}{2}\right\rfloor+\cdots+p(n-km_k,k-1)\\
    &=\sum_{m_k=1}^{\left\lfloor\frac{n}{k}\right\rfloor}p(n-km_k,1)+p(n-km_k,2)+p(n-km_k,3)+\cdots+p(n-km_k,k-1)\\
    &=\sum_{m=1}^{\left\lfloor\frac{n}{k}\right\rfloor}\sum_{\nu=1}^{k-1}p(n-km,\nu).
\end{align*}

Moreover, since
\begin{align}\label{eq:almost}
    p(n)=\sum_{k=1}^np(n,k)=p(n,1)+p(n,2)+\sum_{k=3}^np(n,k)=1+\left\lfloor\frac{n}{2}\right\rfloor+\sum_{k=3}^np(n,k),
\end{align}
we can easily check that $p(0)=1$, $p(1)=1$, $p(2)=2$, and since for every $n\geq3$ \eqref{eq:almost} clearly holds, \eqref{eq:sums_partk1} and \eqref{eq:sums_part1} also hold.
\end{proof}

\section{Generalized smallest parts function}\label{s:spt}
Let $n$ be a nonnegative integer number. Let the smallest part in $\lambda_k^n$ be denoted by $\sigma(\lambda_k^n)$ and the number of smallest parts in $\lambda_k^n$ by $\#(\lambda_k^n)$. The smallest parts function $\spt(n)$ (introduced by Andrews in \cite{andrews}) counts the sum of the number of smallest parts in each partition of $n$, that is, 
\begin{align*}
    \spt(n)=\sum_{\lambda_k^n\in\overline{\lambda^n}}\#(\lambda_k^n).
\end{align*}
We will generalize
\begin{align}\label{eq:sptab}
    \spt_{(a,b)}(n)\coloneqq\sum_{\lambda_k^n\in\overline{\lambda^n}}\sigma(\lambda_k^n)^a\#(\lambda_k^n)^b.
\end{align}

\begin{remark}
    Note that $\spt(n)\coloneqq\spt_{(0,1)}(n)$.
\end{remark}

For example, let $n=5$. We have $\overline{\lambda^5}=\{\{\underbar5\}$, $\{\underbar1,4\}$, $\{\underbar2,3\}$, $\{\underbar1,\underbar1,3\}$, $\{\underbar1,2,2\}$, $\{\underbar1,\underbar1,\underbar1,2\}$, $\{\underbar1,\underbar1,\underbar1,\underbar1,\underbar1\}\}$, where $\sigma(\lambda_k^5)$ is underlined for every partition. $\sigma(\lambda_k^5)$ is, respectively, 5, 1, 2, 1, 1, 1, 1, and $\#(\lambda_k^5)$ is, respectively, 1, 1, 1, 2, 1, 3, 5. Thereby, for instance, \[\spt_{(3,2)}(5)=(5)^3(1)^2+(1)^3(1)^2+(2)^3(1)^2+(1)^3(2)^2+(1)^3(1)^2+(1)^3(3)^2+(1)^3(5)^2=173.\]

We devote this section to developing a recursive and a closed-form formula for \eqref{eq:sptab}.

\begin{lemma}
    \begin{align}\label{eq:spt01}
        \spt_{(1,0)}(n,k)=\sum_{m=1}^{\left\lfloor\frac{n}{k}\right\rfloor}m\sum_{\nu=1}^{k-1}p(n-km,\nu).
    \end{align}
\end{lemma}

\begin{proof}
    Let $\sigma(c_1:\cdots:c_{k-1})_k$ be the sum of the smallest parts of the partitions of $\cc{c_1:\cdots:c_{k-1}}_k$-group. By convention, the smallest part in $\lambda_k^n$ is $\lambda_1$, that is, $\sigma(\lambda_k^n)=\lambda_1=c_1+1$, for which $\sigma(c_1:\cdots:c_{k-1})_k=(c_1+1)\#\cc{c_1:\cdots:c_{k-1}}_k$. From \eqref{eq:c1k},
    \begin{align*}
        \sigma(c_1)_k&=(c_1+1)\#\cc{c_1}_k\\
        &=(c_1+1)\sum_{\nu=1}^{k-1}p(n-k(c_1+1),\nu).
    \end{align*}
    Lastly,
    \begin{align*}
        \spt_{(1,0)}(n,k)&=\sum_{c_1=0}^{\left\lfloor\frac{n}{k}\right\rfloor-1}\sigma(c_1)_k=\sum_{c_1=0}^{\left\lfloor\frac{n}{k}\right\rfloor-1}(c_1+1)\sum_{\nu=1}^{k-1}p(n-k(c_1+1),\nu)=\sum_{m=1}^{\left\lfloor\frac{n}{k}\right\rfloor}m\sum_{\nu=1}^{k-1}p(n-km,\nu).
    \end{align*}
\end{proof}

\begin{corollary}\label{c:spta0}
    Since the $m$ that multiplies the second sum in the right-hand side of \eqref{eq:spt01} is the smallest part of the partitions in $\overline{\lambda_k^n}$, it follows
    \begin{align*}
        \spt_{(a,0)}(n,k)&=\sum_{m=1}^{\left\lfloor\frac{n}{k}\right\rfloor}m^a\sum_{\nu=1}^{k-1}p(n-km,\nu).
    \end{align*}
\end{corollary}

\begin{lemma}
    \begin{align}\label{eq:spt10}
        \spt_{(0,1)}(n,k)=\left\lfloor\frac{k\left\lfloor\frac{n}{k}\right\rfloor}{n}\right\rfloor+\sum_{m=1}^{\left\lfloor\frac{n}{k}\right\rfloor}\sum_{\nu=1}^{k-1}(k-\nu)p(n-km,\nu).
    \end{align}
\end{lemma}

\begin{proof}
    Let $\#[c_1:\cdots:c_{k-1}]_k$ count the number of smallest parts in the partitions of $\cc{c_1:\cdots:c_{k-1}}_k$-group. Let $c_1=c_2=\cdots=c_{k-j-1}=0$ and $1\leq c_{k-j}\leq\cdots\leq c_k$ for $0\leq j\leq k-2$. From \eqref{eq:xparts}, we have $k-j-1$ parts which are equal to 1 (and the rest are greater). Thus, we have
    \begin{align*}
        \#[0:\cdots:0:c_{k-j}:\cdots:c_{k-1}]_k&=(k-j-1)\#\cc{0:\cdots:0:c_{k-j}:\cdots:c_{k-1}}_k\\
        &=(k-j-1)p(n-k,j+1).
    \end{align*}
    Moreover, from \eqref{eq:base1(ii)} we have
    \begin{align*}
        \#[0]_k&=\sum_{j=0}^{k-2}(k-j-1)p(n-k,j+1)\\
        &=\sum_{\nu=1}^{k-1}(k-\nu)p(n-k,\nu).
    \end{align*}
    Since the number of smallest parts does not depend on $c_1$ and it is the same as with $c_1=0$, using \eqref{eq:c1k}, we have
    \begin{align*}
        \#[c_1]_k&=\sum_{\nu=1}^{k-1}(k-\nu)p(n-k(c_1+1),\nu),
    \end{align*}
    and, thus,
    \begin{align}\label{eq:almostspt10}
        \sum_{c_1=0}^{\left\lfloor\frac{n}{k}\right\rfloor-1}\#[c_1]_k=\sum_{c_1=0}^{\left\lfloor\frac{n}{k}\right\rfloor-1}\sum_{\nu=1}^{k-1}(k-\nu)p(n-k(c_1+1),\nu)=\sum_{m=1}^{\left\lfloor\frac{n}{k}\right\rfloor}\sum_{\nu=1}^{k-1}(k-\nu)p(n-km,\nu).
    \end{align}

    If $k\mid n$, there will be one partition $\lambda_k^n=\{\frac nk,\frac nk,\dots,\frac nk\}$ that has one extra smallest part. For example, since $5\mid 15$, the partition $\lambda_5^{15}=\{3,3,3,3,3\}$ has $\#(\lambda_5^{15})=5$ instead of $\#(\lambda_5^{15})=4$ (which is what we consider above). If $k\nmid n$, then we do not have to worry about this part. For example, partition $\lambda_5^{16}=\{3,3,3,3,4\}$ has $\#(\lambda_5^{16})=4$ since $5\nmid 16$.

    If $k\mid n$, we have $\left\lfloor\frac{n}{k}\right\rfloor=\frac nk$. Otherwise, if $k\nmid n$, we have $\left\lfloor\frac{n}{k}\right\rfloor<\frac nk$. Thus, we can add the term
    \begin{align}\label{eq:semid(n)}
        \left\lfloor\frac{\left\lfloor\frac{n}{k}\right\rfloor}{\frac nk}\right\rfloor=\left\lfloor\frac{k\left\lfloor\frac{n}{k}\right\rfloor}{n}\right\rfloor=\begin{cases}
            1 & k\mid n \\
            0 & k\nmid n
        \end{cases}
    \end{align}
    to \eqref{eq:almostspt10} to finalize this proof.
\end{proof}

To proceed, we need \lref{d(n)}.
\begin{lemma}\label{l:d(n)}
    \begin{align*}
        \sum_{k=1}^{n}\left\lfloor\frac{k\left\lfloor\frac{n}{k}\right\rfloor}{n}\right\rfloor=d(n).
    \end{align*}
\end{lemma}

\begin{proof}
    By definition,
    \begin{align*}
        d(n)\coloneqq\sum_{k\mid n}1.
    \end{align*}
    Note that $k$ can only divide $n$ if and only if $1\leq k\leq n$ (in usual circumstances, negative $k$ can also divide $n$, but, by convention, we are only considering positive $k$). We see that if we sum \eqref{eq:semid(n)} through all the possible values of $k$, we are able to count the number of divisors of $n$, that is,
    \begin{align*}
        \sum_{k}\left\lfloor\frac{k\left\lfloor\frac{n}{k}\right\rfloor}{n}\right\rfloor=\sum_{k=1}^n\left\lfloor\frac{k\left\lfloor\frac{n}{k}\right\rfloor}{n}\right\rfloor=\sum_{k\mid n}1+\sum_{k\nmid n}0=d(n).
    \end{align*}
\end{proof}

\begin{corollary}
    \begin{align*}
        \spt(n)=\spt_{(0,1)}(n)=d(n)+\sum_{k=1}^n\sum_{m=1}^{\left\lfloor\frac{n}{k}\right\rfloor}\sum_{\nu=1}^{k-1}(k-\nu)p(n-km,\nu).
    \end{align*}
\end{corollary}

\begin{lemma}
    \begin{align}\label{eq:spt0b}
        \spt_{(0,b)}(n,k)=\left\lfloor\frac{k\left\lfloor\frac{n}{k}\right\rfloor}{n}\right\rfloor(k^b-(k-1)^b)+\sum_{m=1}^{\left\lfloor\frac{n}{k}\right\rfloor}\sum_{\nu=1}^{k-1}(k-\nu)^bp(n-km,\nu).
    \end{align}
\end{lemma}

\begin{proof}
    The second expression of the right-hand side of \eqref{eq:spt0b} follows from the same principle of Corollary \ref{c:spta0}. To see why the first expression holds, we need an example. As before, take $n=5$. We have 
        \begin{align}\label{eq:example1}
            \spt_{(0,b)}(5)=1^b+1^b+1^b+2^b+1^b+3^b+5^b.
        \end{align}
    If we assume
        \begin{align}\label{eq:wrongspt0b}
            \spt_{(0,b)}(n,k)=\left\lfloor\frac{k\left\lfloor\frac{n}{k}\right\rfloor}{n}\right\rfloor^b+\sum_{m=1}^{\left\lfloor\frac{n}{k}\right\rfloor}\sum_{\nu=1}^{k-1}(k-\nu)^bp(n-km,\nu),
        \end{align}
    we have
    \begin{align}
        \notag\spt_{(0,b)}(5)&=\sum_{k=1}^5\left\lfloor\frac{k\left\lfloor\frac{5}{k}\right\rfloor}{5}\right\rfloor^b+\sum_{m=1}^{\left\lfloor\frac{5}{k}\right\rfloor}\sum_{\nu=1}^{k-1}(k-\nu)^bp(5-km,\nu)\\
        \notag&=1^b+\sum_{m=1}^{2}\sum_{\nu=1}^{1}(2-\nu)^bp(5-2m,\nu)+\sum_{m=1}^{1}\sum_{\nu=1}^{2}(3-\nu)^bp(5-3m,\nu)\\
        \notag&\quad+\sum_{m=1}^{1}\sum_{\nu=1}^{3}(4-\nu)^bp(5-4m,\nu)+\left(1^b+\sum_{m=1}^{1}\sum_{\nu=1}^{4}(5-\nu)^bp(5-5m,\nu)\right)\\
        \notag&=1^b+\left(1^bp(3,1)+1^bp(1,1))\right)+\left(2^bp(2,1)+1^bp(2,2)\right)+\left(3^bp(1,1)+2^bp(1,2)+1^bp(1,3)\right)\\
        \notag&\quad+\left(1^b+4^bp(0,1)+3^bp(0,2)+2^bp(0,3)+1^bp(0,4)\right)\\
        \label{eq:example2}&=1^b+1^b+1^b+2^b+1^b+3^b+\left(1^b+4^b\right).
    \end{align}
    Note that the difference between \eqref{eq:example1} and \eqref{eq:example2} is the last term and the last two terms, respectively. \eqref{eq:spt10} calculates the number of smallest parts of each partition considering that there exist, at most, $k-1$ smallest parts. However, any partition $\lambda_k^n=\{\frac nk,\frac nk,\dots,\frac nk\}$ (with $k\mid n$) has a part that is calculated with \eqref{eq:semid(n)}. Thus, when we raise the smallest parts of this partition to the exponent $b$, we get $1^b+(k-1)^b$ instead of $k^b$. To fix this, we could (1) change the entire equation to stop this from occurring, or (2) add the term $k^b-(1^b+(k-1)^b)$ each time $k\mid n$. The second is easier. From \eqref{eq:wrongspt0b} we have
    \begin{align*}
        \left\lfloor\frac{k\left\lfloor\frac{n}{k}\right\rfloor}{n}\right\rfloor^b+\sum_{m=1}^{\left\lfloor\frac{n}{k}\right\rfloor}\sum_{\nu=1}^{k-1}(k-\nu)^bp(n-km,\nu)=\left\lfloor\frac{k\left\lfloor\frac{n}{k}\right\rfloor}{n}\right\rfloor(1)+\sum_{m=1}^{\left\lfloor\frac{n}{k}\right\rfloor}\sum_{\nu=1}^{k-1}(k-\nu)^bp(n-km,\nu)
    \end{align*}
    and adding our term we get
    \begin{align*}
        \spt_{(0,b)}(n,k)&=\left\lfloor\frac{k\left\lfloor\frac{n}{k}\right\rfloor}{n}\right\rfloor\left(1+k^b-(1^b+(k-1)^b)\right)+\sum_{m=1}^{\left\lfloor\frac{n}{k}\right\rfloor}\sum_{\nu=1}^{k-1}(k-\nu)^bp(n-km,\nu)\\
        &=\left\lfloor\frac{k\left\lfloor\frac{n}{k}\right\rfloor}{n}\right\rfloor\left(k^b-(k-1)^b)\right)+\sum_{m=1}^{\left\lfloor\frac{n}{k}\right\rfloor}\sum_{\nu=1}^{k-1}(k-\nu)^bp(n-km,\nu).
    \end{align*}
\end{proof}

\begin{lemma}
    \begin{align}\label{eq:spt11}
        \spt_{(1,1)}(n,k)=\left\lfloor\frac{n}{k}\right\rfloor\left\lfloor\frac{k\left\lfloor\frac{n}{k}\right\rfloor}{n}\right\rfloor+\sum_{m=1}^{\left\lfloor\frac{n}{k}\right\rfloor}m\sum_{\nu=1}^{k-1}(k-\nu)p(n-km,\nu).
    \end{align}
\end{lemma}

\begin{proof}
    When $k\mid n$, each part of partition $\lambda_k^n=\{\frac nk,\frac nk,\dots,\frac nk\}$ is equal to $\frac nk$. Since the first expression in \eqref{eq:spt11} counts these partitions, we may multiply $\frac nk$ to get the product $\sigma(\lambda_k^n)\#(\lambda_k^n)$ for this partition. The second term follows immediately, since \eqref{eq:sptab} is multiplicative. 
\end{proof}

\begin{corollary}
    \begin{align}
        \label{eq:sptabnk}\spt_{(a,b)}(n,k)&=\left\lfloor\frac{n}{k}\right\rfloor^a\left\lfloor\frac{k\left\lfloor\frac{n}{k}\right\rfloor}{n}\right\rfloor\left(k^b-(k-1)^b)\right)+\sum_{m=1}^{\left\lfloor\frac{n}{k}\right\rfloor}m^a\sum_{\nu=1}^{k-1}(k-\nu)^bp(n-km,\nu),\\
        \notag\spt_{(a,b)}(n)&=\sum_{k=1}^n\left\lfloor\frac{n}{k}\right\rfloor^a\left\lfloor\frac{k\left\lfloor\frac{n}{k}\right\rfloor}{n}\right\rfloor\left(k^b-(k-1)^b)\right)+\sum_{m=1}^{\left\lfloor\frac{n}{k}\right\rfloor}m^a\sum_{\nu=1}^{k-1}(k-\nu)^bp(n-km,\nu).
    \end{align}
\end{corollary}

\begin{remark}
    We have written all formulas related to the smallest parts of partitions in recursive form. However, one may use \eqref{eq:sums_partk1} in any of these formulas to make them closed-form.
\end{remark}

\begin{theorem}\label{th:ineq}
    $p(n,k)\leq\spt_{(a,b)}(n,k)$ and $p(n)\leq\spt_{(a,b)}(n)$ for $a,b\geq0$.
\end{theorem}

\begin{proof}
    Firstly, from \eqref{eq:rec_partk1} and \eqref{eq:sptabnk}, we have
    \begin{align*}
        p(n,k)&\leq\spt_{(a,b)}(n,k)\\
        \sum_{m=1}^{\left\lfloor\frac{n}{k}\right\rfloor}\sum_{\nu=1}^{k-1}p(n-km,\nu)&\leq\left\lfloor\frac{n}{k}\right\rfloor^a\left\lfloor\frac{k\left\lfloor\frac{n}{k}\right\rfloor}{n}\right\rfloor\left(k^b-(k-1)^b)\right)+\sum_{m=1}^{\left\lfloor\frac{n}{k}\right\rfloor}m^a\sum_{\nu=1}^{k-1}(k-\nu)^bp(n-km,\nu).
    \end{align*}
    If we let $a=b=0$, we have
    \begin{align*}
        \sum_{m=1}^{\left\lfloor\frac{n}{k}\right\rfloor}\sum_{\nu=1}^{k-1}p(n-km,\nu)&=\left\lfloor\frac{n}{k}\right\rfloor^0\left\lfloor\frac{k\left\lfloor\frac{n}{k}\right\rfloor}{n}\right\rfloor\left(k^0-(k-1)^0)\right)+\sum_{m=1}^{\left\lfloor\frac{n}{k}\right\rfloor}m^0\sum_{\nu=1}^{k-1}(k-\nu)^0p(n-km,\nu)\\
        \sum_{m=1}^{\left\lfloor\frac{n}{k}\right\rfloor}\sum_{\nu=1}^{k-1}p(n-km,\nu)&=\sum_{m=1}^{\left\lfloor\frac{n}{k}\right\rfloor}\sum_{\nu=1}^{k-1}p(n-km,\nu)\\
        p(n,k)&=\spt_{(a,b)}(n,k).
    \end{align*}
    We use the convention that $p(0,1)=1$, and thus we can also use the convention that $\spt_{(a,b)}(0,1)=1$. Moreover, if we let $n=1$, we have
    \begin{align*}
        p(1,k)&=\spt_{(a,b)}(1,k)\\
        1&=\left\lfloor\frac{1}{1}\right\rfloor^a\left\lfloor\frac{1\left\lfloor\frac{1}{1}\right\rfloor}{1}\right\rfloor\left(1^b-(1-1)^b)\right)+\sum_{m=1}^{\left\lfloor\frac{1}{1}\right\rfloor}m^a\sum_{\nu=1}^{1-1}(1-\nu)^bp(1-m,\nu)\\
        1&=1.
    \end{align*}
    Additionally, if we let $a,b>0$ and $n>1$, note that \[\left\lfloor\frac{n}{k}\right\rfloor^a\left\lfloor\frac{k\left\lfloor\frac{n}{k}\right\rfloor}{n}\right\rfloor\left(k^b-(k-1)^b)\right),m^a,(k-\nu)^b>1,\] for we clearly see that $p(n,k)<\spt_{(a,b)}(n,k)$.
    
    Lastly, we have
    \begin{align*}
        p(n)&\leq\spt_{(a,b)}(n)\\
        \sum_{k=1}^np(n,k)&\leq\sum_{k=1}^n\spt_{(a,b)}(n,k)\\
        p(n,k)&\leq\spt_{(a,b)}(n,k)
    \end{align*}
    which we have seen holds.
\end{proof}

\section{Relatively and not-relatively prime partition function}\label{s:lambda}
In this section, we give an intuitive formula for the not-relatively prime partition function by relying on the inclusion-exclusion principle, and then we use it to derive a formula for the relatively prime partition function. 

Let $\Lambda_m(n)$ denote the number of partitions of $n$ in which all parts are divisible by $m$. Let $\overline{\Lambda_m(n)}$ be the set of all partitions in which all parts are divisible by $m$ and let $\overline{\Lambda_m(n,k)}$ be the equivalent for $\Lambda_m(n,k)$. Note that $\left|\overline{\Lambda_m(n)}\right|=\Lambda_m(n)$ and $\left|\overline{\Lambda_m(n,k)}\right|=\Lambda_m(n,k)$. 

\begin{theorem}
    The following identities hold:
    \begin{align*}
        \Lambda_m(n,k)=\begin{cases}
            p\left(\frac nm,k\right) & m\mid n \\
            0 & m\nmid n
        \end{cases}\qquad\text{and}\qquad\Lambda_m(n)=\begin{cases}
            p\left(\frac nm\right) & m\mid n \\
            0 & m\nmid n
        \end{cases}.
    \end{align*}
\end{theorem}

\begin{proof}
    Due to Corollary \ref{c:general}, if $m\nmid n$ then it directly follows that there are no partitions of $n$ in which all parts are divisible by $m$, and thus $\Lambda_m(n,k)=\Lambda_m(n)=0$. Otherwise, if $m\mid n$, then there exists some integer $a$ such that
    \begin{align*}
        n&=am\\
        \frac nm&=a\\
        p\left(\frac nm,k\right)&=p(a,k).
    \end{align*}
    
    Since $p(a,k)$ counts all the partitions of $n$ in which all parts are divisible by $m$, it follows that $\Lambda_m(n,k)=p\left(\frac nm,k\right)$, and thus $\Lambda_m(n)=p\left(\frac nm\right)$.
\end{proof}

Recall that the inclusion-exclusion principle states that the cardinality of the union of a finite amount of finite sets is the sum of the cardinality of intersections with an odd amount of sets minus the sum of the cardinality of intersections with an even amount of sets. That is, for finite sets $A_1,\dots,A_\ell$,
\begin{align}\label{eq:inex}
    \left|\bigcup_{\nu=1}^\ell A_\nu\right|&=\sum_{i=1}^\ell(-1)^{i+1}\sum_{1\leq \nu_1<\cdots<\nu_i\leq \ell}\left|\bigcap_{j=1}^iA_{\nu_j}\right|\\
    \notag&=\sum_{1\leq\nu_1\leq\ell}|A_{\nu_1}|-\sum _{1\leq\nu_1<\nu_2\leq\ell}|A_{\nu_1}\cap A_{\nu_2}|+
    \cdots+(-1)^{\ell+1}\sum _{1\leq\nu_1<\cdots<\nu_\ell\leq\ell}\left|A_{\nu_1}\cap\cdots\cap A_{\nu_\ell}\right|.
\end{align}

Due to Corollary \ref{c:general}, the only not-relatively prime partitions of some integer $n$ are those in which all parts are divisible by some integer $1\ne m\mid n$, and thus it follows that
\begin{align*}
    \Lambda(n,k)&=\left|\bigcup_{1\neq m\mid n}\overline{\Lambda_m(n,k)}\right|=\sum_{i=1}^{d(n)-1}(-1)^{i+1}\sum_{\substack{1<m_1<\cdots<m_{i}\leq n \\ m_1,\dots,m_{i}\mid n}}\left|\bigcap_{j=1}^i\overline{\Lambda_{m_j}(n,k)}\right|.
\end{align*}
We exclude $m=1$, even though 1 is a factor of all integers, since $\overline{\Lambda_1(n,k)}$ contains relatively prime partitions, and thus there are some partitions in $\overline{\Lambda_1(n,k)}$ that are not in $\overline{\Lambda(n,k)}$. (In fact, note that $\left|\overline{\Lambda_1(n,k)}\right|=\Lambda_1(n,k)=p\left(\frac n1,k\right)=p(n,k)\geq\Lambda(n,k)$). Note that this formula is not defined for $n=1$ since the only factor of 1 is 1, and we exclude it. However, $\Lambda(1,k)$ and thus $\Lambda(1)$ is defined in \eqref{eq:base_relprime}. $d(n)-1$ appears here since $n$ has $d(n)-1$ divisors without considering 1, and thus analogous to \eqref{eq:inex}, $\ell=d(n)-1$. 

The intersection $\overline{\Lambda_{m_1}(n,k)}\cap\overline{\Lambda_{m_2}(n,k)}\cap\cdots\cap\overline{\Lambda_{m_i}(n,k)}=\overline{\Lambda_{m'}(n,k)}$ such that $m'$ is the smallest number that is divisible by each of $m_1,\dots,m_i$. Note that, by definition, $m'=\lcm(m_1,\dots,m_i)$, and thus
\begin{align*}
    \left|\bigcap_{j=1}^i\overline{\Lambda_{m_j}(n,k)}\right|&=\left|\overline{\Lambda_{\lcm(m_1,\dots,m_i)}(n,k)}\right|=\Lambda_{\lcm(m_1,\dots,m_i)}(n,k)=p\left(\frac n{\lcm(m_1,\dots,m_i)},k\right).
\end{align*}
As an example, note that 
\begin{align*}
    \overline{\Lambda_{12}(420,2)}\cap\overline{\Lambda_{42}(420,2)}=\overline{\Lambda_{\lcm(12,42)}(420,2)}=\overline{\Lambda_{84}(420,2)}=\{\{84,336\},\{168,252\}\}
\end{align*}
and hence
\begin{align*}
    \left|\overline{\Lambda_{12}(420,2)}\cap\overline{\Lambda_{42}(420,2)}\right|=\left|\overline{\Lambda_{\lcm(12,42)}(420,2)}\right|=\left|\overline{\Lambda_{84}(420,2)}\right|=p\left(\frac {420}{84},2\right)=\left\lfloor\frac{\frac {420}{84}}2\right\rfloor=2.
\end{align*}

Since $p(n,2)=\left\lfloor\frac{n}{2}\right\rfloor$, if $n<2$, then it follows that $p(n,2)=0$. Otherwise, for $k\geq3$, $p(n,k)=\sum_{m_k=1}^{\left\lfloor\frac{n}{k}\right\rfloor}\cdots$, and thus if $\left\lfloor\frac{n}{k}\right\rfloor<1$, then $p(n,k)=0$. Note that $\left\lfloor\frac{a}{b}\right\rfloor<1$ if $a<b$ (for nonnegative $a$ and positive $b$), for which we require
\begin{align*}
    \frac n{\lcm(m_1,\dots,m_i)}&\geq k\\
    \lcm(m_1,\dots,m_i)&\leq\frac nk.
\end{align*}

Putting everything together, we have
\begin{align*}
    \Lambda(n,k)&=\sum_{i\geq1}(-1)^{i+1}\sum_{\substack{1<m_1<\cdots<m_{i}\leq n \\ m_1,\dots,m_{i}\mid n \\ \lcm(m_1,\dots,m_i)\leq\frac nk}}p\left(\frac n{\lcm(m_1,\dots,m_i)},k\right),
\end{align*}
which can easily be extended to
\begin{align*}
    \Lambda(n)&=\sum_{k=1}^n\Lambda(n,k)=\sum_{k=1}^n\sum_{i\geq1}(-1)^{i+1}\sum_{\substack{1<m_1<\cdots<m_{i}\leq n \\ m_1,\dots,m_{i}\mid n \\ \lcm(m_1,\dots,m_i)\leq\frac nk}}p\left(\frac n{\lcm(m_1,\dots,m_i)},k\right).
\end{align*}

\begin{remark}
    These expressions are neither efficient nor in closed-form. However, they allow us to derive a formula for the relatively prime partition function. 
\end{remark}

Due to Corollary \ref{c:formula}, we have $p_\Psi(n,k)=p(n,k)-\Lambda(n,k)$ and $p_\Psi(n)=p(n)-\Lambda(n)$, for we see that
\begin{align}
    &\notag p_\Psi(n,2)=\left\lfloor\frac{n}{2}\right\rfloor+\sum_{i\geq1}(-1)^{i}\sum_{\substack{1<m_1<\cdots<m_{i}\leq n \\ m_1,\dots,m_{i}\mid n \\ \lcm(m_1,\dots,m_i)\leq\frac nk}}p\left(\frac n{\lcm(m_1,\dots,m_i)},k\right),\\
    \notag &p_\Psi(n,k)=\sum_{m=1}^{\left\lfloor\frac{n}k\right\rfloor}\sum_{\nu=1}^{k-1}p(n-km,\nu)+\sum_{i\geq1}(-1)^{i}\sum_{\substack{1<m_1<\cdots<m_{i}\leq n \\ m_1,\dots,m_{i}\mid n \\ \lcm(m_1,\dots,m_i)\leq\frac nk}}p\left(\frac n{\lcm(m_1,\dots,m_i)},k\right)\\
    \label{eq:relprimek}&=\sums_{\substack{t=0 \\ \left[k-t,\,1,\,0,\,\left\lfloor\frac{n-\sum_{j=k-t+1}^{k}jm_j}{k-t}\right\rfloor\right]}}^{k-3}\left\lfloor\frac{2+n-\sum_{j=3}^{k}jm_j}{2}\right\rfloor+\sum_{i\geq1}(-1)^{i}\sum_{\substack{1<m_1<\cdots<m_{i}\leq n \\ m_1,\dots,m_{i}\mid n \\ \lcm(m_1,\dots,m_i)\leq\frac nk}}p\left(\frac n{\lcm(m_1,\dots,m_i)},k\right),
\end{align}
for $k\geq3$, and, for $n>1$,
\begin{align*}
    &=\sum_{m=1}^{\left\lfloor\frac{n}k\right\rfloor}\sum_{\nu=1}^{k-1}p(n-km,\nu)+\sum_{k=1}^n\sum_{i\geq1}(-1)^{i}\sum_{\substack{1<m_1<\cdots<m_{i}\leq n \\ m_1,\dots,m_{i}\mid n \\ \lcm(m_1,\dots,m_i)\leq\frac nk}}p\left(\frac n{\lcm(m_1,\dots,m_i)},k\right)\\
    &=1+\left\lfloor\frac n2\right\rfloor+\sum_{k=3}^n\sums_{\substack{t=0 \\ \left[k-t,\,1,\,0,\,\left\lfloor\frac{n-\sum_{j=k-t+1}^{k}jm_j}{k-t}\right\rfloor\right]}}^{k-3}\left\lfloor\frac{2+n-\sum_{j=3}^{k}jm_j}{2}\right\rfloor
    \\\notag&\quad+\sum_{k=1}^n\sum_{i\geq1}(-1)^{i}\sum_{\substack{1<m_1<\cdots<m_{i}\leq n \\ m_1,\dots,m_{i}\mid n \\ \lcm(m_1,\dots,m_i)\leq\frac nk}}p\left(\frac n{\lcm(m_1,\dots,m_i)},k\right).
\end{align*}

As an example, we may look at the relatively and not-relatively prime partitions of 30 with 2 parts. We see that the factors of 30 (without considering 1) are 2, 3, 5, 6, 10, 15, and 30. Thus, we have
\begin{align*}
    \Lambda(30,2)&=\sum_{i\geq1}(-1)^{i+1}\sum_{\substack{1<m_1<\cdots<m_{i}\leq 30 \\ m_1,\dots,m_{i}\mid 30 \\ \lcm(m_1,\dots,m_i)\leq\frac {30}2}}p\left(\frac {30}{\lcm(m_1,\dots,m_i)},2\right)\\
    &=\left(p\left(\frac{30}{2},2\right)+p\left(\frac{30}{3},2\right)+p\left(\frac{30}{5},2\right)+p\left(\frac{30}{6},2\right)+p\left(\frac{30}{10},2\right)+p\left(\frac{30}{15},2\right)\right)
    \\&\quad-\left(p\left(\frac{30}{\lcm(2,3)},2\right)+p\left(\frac{30}{\lcm(2,5)},2\right)+p\left(\frac{30}{\lcm(2,6)},2\right)+p\left(\frac{30}{\lcm(2,10)},2\right)\right.
    \\&\quad+\left.p\left(\frac{30}{\lcm(3,5)},2\right)+p\left(\frac{30}{\lcm(3,6)},2\right)+p\left(\frac{30}{\lcm(3,15)},2\right)+p\left(\frac{30}{\lcm(5,10)},2\right)\right.
    \\&\quad+\left.p\left(\frac{30}{\lcm(5,15)},2\right)\right)+\left(p\left(\frac{30}{\lcm(2,3,6)},2\right)+p\left(\frac{30}{\lcm(2,5,10)},2\right)+p\left(\frac{30}{\lcm(3,5,15)},2\right)\right)\\
    &=\left(\left\lfloor\frac{15}{2}\right\rfloor+\left\lfloor\frac{10}{2}\right\rfloor+\left\lfloor\frac{6}{2}\right\rfloor+\left\lfloor\frac{5}{2}\right\rfloor+\left\lfloor\frac{3}{2}\right\rfloor+\left\lfloor\frac{2}{2}\right\rfloor\right)-\left(\left\lfloor\frac{5}{2}\right\rfloor+\left\lfloor\frac{3}{2}\right\rfloor+\left\lfloor\frac{5}{2}\right\rfloor+\left\lfloor\frac{3}{2}\right\rfloor+\left\lfloor\frac{2}{2}\right\rfloor\right.
    \\&\quad+\left.\left\lfloor\frac{5}{2}\right\rfloor+\left\lfloor\frac{2}{2}\right\rfloor+\left\lfloor\frac{3}{2}\right\rfloor+\left\lfloor\frac{2}{2}\right\rfloor\right)+\left(\left\lfloor\frac{5}{2}\right\rfloor+\left\lfloor\frac{3}{2}\right\rfloor+\left\lfloor\frac{2}{2}\right\rfloor\right)\\
    &=(7+5+3+2+1+1)-(2+1+2+1+1+2+1+1+1)+(2+1+1)=11
\end{align*}
and
\begin{align*}
    p_\Psi(30,2)=p(30,2)-\Lambda(30,2)=\left\lfloor\frac{30}{2}\right\rfloor-11=4.
\end{align*}

\section{Links to totient functions}\label{s:imp}
We can extend \eqref{eq:relprimek} to derive identities for Jordan's totient function of second order, Euler's totient function, and Dedekind's psi function. In 2008, El Bachraoui \cite{elbachraoui} proved the following identities:
\begin{align}
    \label{eq:bachraoui1}p_\Psi(n,2)&=\frac n2\prod_{p\mid n}\left(1-\frac1p\right)=\frac12\phi(n),\\
    \label{eq:bachraoui2}p_\Psi(n,3)&=\frac{n^2}{12}\prod_{p\mid n}\left(1-\frac1{p^2}\right)=\frac1{12}J_2(n).
\end{align}
We will link these identities with the partition function and the not-relatively prime partition function, alike, and develop identities for $J_2(n)$, $\phi(n)$, and $\psi(n)$. Before starting, recall that
\begin{equation}\label{eq:identities}\begin{tabular}{c}
    $\displaystyle\Lambda(n,1)=\begin{cases} 0 & n=1 \\ 1 & n\neq1 \end{cases};\quad p_\Psi(n,1)=\begin{cases} 1 & n=1 \\ 0 & n\neq1 \end{cases};\quad p(n,k)=p_\Psi(n,k)+\Lambda(n,k);$\\
    $\displaystyle p(n,1)=1;\quad p(n,2)=\left\lfloor\frac n2\right\rfloor;\quad p(n,3)=\sum_{m_3=1}^{\left\lfloor\frac n3\right\rfloor}\left\lfloor\frac{2+n-3m_3}2\right\rfloor;$\\
    $\displaystyle\Lambda(n,2)=\sum_{i\geq1}(-1)^{i+1}\sum_{\substack{1<m_1<\cdots<m_{i}\leq n \\ m_1,\dots,m_{i}\mid n \\ \lcm(m_1,\dots,m_i)\leq\frac n2}}p\left(\frac n{\lcm(m_1,\dots,m_i)},2\right);$\\
    $\displaystyle\Lambda(n,3)=\sum_{i\geq1}(-1)^{i+1}\sum_{\substack{1<m_1<\cdots<m_{i}\leq n \\ m_1,\dots,m_{i}\mid n \\ \lcm(m_1,\dots,m_i)\leq\frac n3}}p\left(\frac n{\lcm(m_1,\dots,m_i)},3\right)$.
\end{tabular}\end{equation}

Going forward, we will only write $p(n,1)$, $p(n,2)$, $p(n,3)$, $\Lambda(n,1)$, $\Lambda(n,2)$, and $p_\Psi(n,1)$, and we will refer the reader to \eqref{eq:identities}. 

From \eqref{eq:bachraoui1} and \eqref{eq:bachraoui2}, we respectively have
\begin{align}
    \notag p_\Psi(n,2)&=p(n,2)-\Lambda(n,2)\\
    \label{eq:elbachraoui1}\phi(n)&=2(p(n,2)-\Lambda(n,2)),\\
    \notag p_\Psi(n,3)&=p(n,3)-\Lambda(n,3)\\
    \label{eq:elbachraoui2}J_2(n)&=12(p(n,3)-\Lambda(n,3)).
\end{align}
In 1920, Hardy \cite{hardy} showed that $p(n,1)+p(n,2)+p(n,3)=\left\langle\frac{(n+3)^2}{12}\right\rangle$. We have
\begin{align*}
    p_\Psi(n,1)&=p(n,1)-\Lambda(n,1),\\
    p_\Psi(n,2)&=p(n,2)-\Lambda(n,2),\\
    p_\Psi(n,3)&=p(n,3)-\Lambda(n,3),
\end{align*}
and, thus,
\begin{align*}
    p_\Psi(n,1)+p_\Psi(n,2)&=p(n,1)-\Lambda(n,1)+p_\Psi(n,2)\\
    p_\Psi(n,1)+p_\Psi(n,2)&=p(n,1)+p(n,2)-\Lambda(n,1)-\Lambda(n,2)\\
    p_\Psi(n,1)+p_\Psi(n,2)+p_\Psi(n,3)&=p(n,1)+p(n,2)-\Lambda(n,1)-\Lambda(n,2)+p_\Psi(n,3)\\
    p_\Psi(n,1)+p_\Psi(n,2)+p_\Psi(n,3)&=p(n,1)+p(n,2)+p(n,3)-\Lambda(n,1)-\Lambda(n,2)-\Lambda(n,3),
\end{align*}
for we see that, if $n=1$, we have
\begin{align*}
    1+\frac12\phi(n)+\frac1{12}J_2(n)&=\left\langle\frac{(n+3)^2}{12}\right\rangle-\Lambda(n,2)-\Lambda(n,3),
\end{align*}
and, if $n\neq1$, we have
\begin{align*}
    \frac12\phi(n)+\frac1{12}J_2(n)&=\left\langle\frac{(n+3)^2}{12}\right\rangle-1-\Lambda(n,2)-\Lambda(n,3).
\end{align*}
Note that both expressions are equal. From this, we have
\begin{align}
    \label{eq:hardy1}\phi(n)&=2\left(\left\langle\frac{(n+3)^2}{12}\right\rangle-\Lambda(n,2)-\Lambda(n,3)-1-\frac1{12}J_2(n)\right)\\
    \label{eq:hardy2}J_2(n)&=12\left(\left\langle\frac{(n+3)^2}{12}\right\rangle-\Lambda(n,2)-\Lambda(n,3)-1-\frac12\phi(n)\right).
\end{align}

Furthermore, in 1985, Honsberger \cite{honsberger}\footnote{This same result can also be found in \cite[eqs. (67) and (69)]{weisstein}.} showed that $p(n,3)=\left\langle\frac{n^2}{12}\right\rangle$, for we can see that, following from \eqref{eq:elbachraoui2}, we have
\begin{align*}
    J_2(n)&=12\left(\left\langle\frac{n^2}{12}\right\rangle-\Lambda(n,3)\right).
\end{align*}

In \cite[eqs. (62), (63), (67), and (68)]{weisstein} the following identities are covered:
\begin{align*}
    p(n,2)&=\left\langle\frac{2n-1}{4}\right\rangle,\\
    p(n,2)&=\frac14\left\langle2n-1+(-1)^n\right\rangle\\
    &=\frac{2n-1+(-1)^n}4,\\
    p(n,3)&=\frac1{72}\left\langle6n^2-7-9(-1)^n+16\cos\left(\frac{2\pi n}3\right)\right\rangle\\
    &=\frac{6n^2-7-9(-1)^n+16\cos\left(\frac{2\pi n}3\right)}{72},
\end{align*}
for which, using \eqref{eq:elbachraoui1} and \eqref{eq:elbachraoui2}, we respectively get
\begin{equation*}\begin{split}
    \phi(n)&=2\left(\left\langle\frac{2n-1}{4}\right\rangle-\Lambda(n,2)\right),\\
    \phi(n)&=2\left(\frac14\left\langle2n-1+(-1)^n\right\rangle-\Lambda(n,2)\right),\\
    &=n+\frac{(-1)^n-1}2-2\Lambda(n,2)\\
    J_2(n)&=12\left(\frac1{72}\left\langle6n^2-7-9(-1)^n+16\cos\left(\frac{2\pi n}3\right)\right\rangle-\Lambda(n,3)\right)\\
    &=n^2+\frac{16\cos\left(\frac{2\pi n}3\right)-7-9(-1)^n}{6}-12\Lambda(n,3).
\end{split}\end{equation*}

We can combine \eqref{eq:elbachraoui2} and \eqref{eq:elbachraoui1} with \eqref{eq:hardy1} and \eqref{eq:hardy2}, respectively, to get
\begin{align*}
    \phi(n)&=2\left(\left\langle\frac{(n+3)^2}{12}\right\rangle-\Lambda(n,2)-\Lambda(n,3)-1-\frac1{12}(12(p(n,3)-\Lambda(n,3)))\right)\\
    &=2\left(\left\langle\frac{(n+3)^2}{12}\right\rangle-\Lambda(n,2)-\Lambda(n,3)-1-p(n,3)+\Lambda(n,3)\right)\\
    &=2\left(\left\langle\frac{(n+3)^2}{12}\right\rangle-\Lambda(n,2)-1-p(n,3)\right),\\
    J_2(n)&=12\left(\left\langle\frac{(n+3)^2}{12}\right\rangle-\Lambda(n,2)-\Lambda(n,3)-1-\frac12(2(p(n,2)-\Lambda(n,2)))\right)\\
    &=12\left(\left\langle\frac{(n+3)^2}{12}\right\rangle-\Lambda(n,2)-\Lambda(n,3)-1-p(n,2)+\Lambda(n,2)\right)\\
    &=12\left(\left\langle\frac{(n+3)^2}{12}\right\rangle-\Lambda(n,3)-1-p(n,2)\right).
\end{align*}
With this and the results from \cite{honsberger,weisstein}, we get
\begin{align*}
    \phi(n)&=2\left(\left\langle\frac{(n+3)^2}{12}\right\rangle-\Lambda(n,2)-1-\left\langle\frac{n^2}{12}\right\rangle\right),\\
    \phi(n)&=2\left(\left\langle\frac{(n+3)^2}{12}\right\rangle-\Lambda(n,2)-1-\frac{6n^2-7-9(-1)^n+16\cos\left(\frac{2\pi n}3\right)}{72}\right),\\
    J_2(n)&=12\left(\left\langle\frac{(n+3)^2}{12}\right\rangle-\Lambda(n,3)-1-\left\langle\frac{2n-1}{4}\right\rangle\right),\\
    J_2(n)&=12\left(\left\langle\frac{(n+3)^2}{12}\right\rangle-\Lambda(n,3)-1-\frac{2n-1+(-1)^n}4\right)\\
    &=3\left(4\left\langle\frac{(n+3)^2}{12}\right\rangle-4\Lambda(n,3)-5-2n+(-1)^n\right).
\end{align*}

Note that
\begin{align*}
    \psi(n)=\frac{J_2(n)}{\phi(n)}&=\frac22\left(\frac{J_2(n)}{\phi(n)}\right)=\frac12\left(\frac{J_2(n)}{p_\Psi(n,2)}\right)=\frac12\left(\frac{J_2(n)}{p(n,2)-\Lambda(n,2)}\right)\\
    &=\frac{12}{12}\left(\frac{J_2(n)}{\phi(n)}\right)=12\left(\frac{p_\Psi(n,3)}{\phi(n)}\right)=12\left(\frac{p(n,3)-\Lambda(n,3)}{\phi(n)}\right).
\end{align*}
We can combine all the different identities of $\phi(n)$ and $J_2(n)$ that we have discussed in this section to develop numerous identities for $\psi(n)$. Furthermore, we can also create new identities involving $\phi(n)$, $J_2(n)$, and $\psi(n)$ for $p(n,2)$, $p(n,3)$, $p_\Psi(n,2)$, $p_\Psi(n,3)$, $\Lambda(n,2)$, and $\Lambda(n,3)$ by solving for them in the various identities we have discovered. 
Other relations to the partition function can similarly lay new identities for functions related to the partitions of an integer. 

\section{Further work}
We remark the importance of \thref{part} for various problems involving the enumeration of restricted partitions, as well as the study of individual partitions themselves. The notions of base and applied representations of $k$ in \sref{formula} can be extensively used to develop closed-form formulas for most functions of restricted partitions (such as $p_\Psi(n)$, $\Lambda(n)$, etc.) analogous to how the generating function of $p(n)$ (see \eqref{eq:euler} below) is commonly used to derive new generating functions for functions of restricted partitions. 

It is an interesting problem to ask for which values of $a,a',b,b'$
\begin{align*}
    \spt_{(a,b)}(n)\leq\spt_{(a',b')}(n)
\end{align*}
holds and for which values of $a,a',b,b',C,k,k'$
\begin{align*}
    \spt_{(a,b)}(n,k)\leq C\spt_{(a,b)}(n,k');\quad\spt_{(a,b)}(n,k)\leq C\spt_{(a',b')}(n,k')
\end{align*}
hold. For instance, since $\lambda_1^n=\{n\}$ and $\lambda_n^n=\{1,1,\dots,1\}$, we can see that 
\begin{align*}
    \spt_{(a,0)}(n,1)&=\spt_{(0,a)}(n,n)=n^a\spt_{(0,1)}(n,1)=n^a\spt_{(1,0)}(n,n)=n^a.
\end{align*}

Recall Euler's famous generating function
\begin{align}\label{eq:euler}
    \sum_{n=1}^\infty p(n)q^n&=\frac1{(q)_\infty},
\end{align}
where
\begin{align*}
    (q)_n\coloneqq\prod_{m=1}^n(1-q^m);\quad(q)_\infty\coloneqq\lim_{n\to\infty}(q)_n;\quad(q)_0=1.
\end{align*}
\cite[Theorems 2.4 and 3.1]{spt} provides the generating functions
\begin{align*}
    \sum_{n=1}^\infty\spt_{(0,1)}(n)q^n&=\frac1{(q)_\infty}\sum_{n=1}^\infty\frac{q^n(q)_{n-1}}{1-q^n},\\
    \sum_{n=1}^\infty\spt_{(1,1)}(n)q^n&=\frac1{(q)_\infty}\sum_{n=1}^\infty\frac{nq^n(q)_{n-1}}{1-q^n}.
\end{align*}

\begin{conjecture}
    We may suspect that
    \begin{align*}
        \sum_{n=1}^\infty\spt_{(a,b)}(n)q^n&=\frac1{(q)_\infty}\sum_{n=1}^\infty n^a\left(\frac{q^n(q)_{n-1}}{1-q^n}\right)^b.
    \end{align*}
\end{conjecture}

If a closed-form formula for $p_\Psi(n)$ is uncovered (using, for instance, the methods of \sref{formula}), then there will be many more new identities for the totient functions studied in \sref{imp}. Moreover, more identities can be found by relating these totient functions to other problems or functions unrelated to partitions. 

Lastly, we offer a generalization of $\xi$-notation. Let $\star$ be some arbitrary operator. We may use circle-notation in the following manner:
\begin{align}\label{eq:circ}
    \circs_{\substack{P(j) \\ \langle\star\rangle}}a_j&=a_1\star a_2\star\cdots\star a_{Q(j)}.
\end{align}
We would say that the left-hand side of \eqref{eq:circ} is a circle operator of $\star$, and we would abbreviate it by writing $\TikCircle[1.3]_\star$. Similarly, if we have some big operator $\square$, we may write
\begin{align*}
    \circs_{\substack{P(j) \\ [parameters] \\ \langle\square\rangle}}C&=\squares_{P_1(\nu)}\squares_{P_2(\nu)}\cdots\squares_{P_{Q(j)}(\nu)}C.
\end{align*}
To avoid confusion, there must be explanation of which operator will be in use.

For example, we may write $\Sigma$-notation as
\begin{align*}
    \circs_{\substack{P(j) \\ \langle+\rangle}}a_j&=a_1+ a_2+\cdots+a_{Q(j)}=\sum_{P(j)}a_j
\end{align*}
or $\xi$-notation as
\begin{align*}
    \circs_{\substack{P(j) \\ [parameters] \\ \langle\Sigma\rangle}}C&=\sum_{P_1(\nu)}\sum_{P_2(\nu)}\cdots\sum_{P_{Q(j)}(\nu)}C=\sums_{\substack{P(j) \\ [parameters] \\ \langle\Sigma\rangle}}C.
\end{align*}

\section*{}

\end{document}